\documentclass[12pt]{article}
\let\proof\relax
 \let\endproof\relax

 \usepackage{amsmath}
\usepackage{graphicx}
\usepackage{caption,subcaption}
\usepackage[utf8]{inputenc}
  \usepackage{nomencl}

\usepackage{placeins} 
\usepackage{flafter}  

\usepackage[misc,geometry]{ifsym} 
 
 \usepackage{xcolor}
\usepackage{ctable}
\usepackage{datetime}
\usepackage{amsthm, amssymb}
\usepackage{lineno}
\usepackage{epstopdf}

\usepackage{graphicx,float}
\usepackage{color}
\usepackage{amsmath}
\usepackage{amssymb}
\usepackage[normalem]{ulem}

\newcommand{\bF}{{\mathbf F}}
\newcommand{\bK}{{\mathbf K}}

\newcommand{\bI}{{\mathbf I}}

\newcommand{\bu}{{\mathbf u}}

\newcommand{\be}{{\mathbf e}}

\newcommand{\bb}{\mathfrak{b}}
\newcommand{\cc}{\mathfrak{c}}
\newcommand{\ii}{i}
\newcommand{\iim}{i-1}

\newcommand{\ga}{\mathbf{g}} 
\newcommand{\uu}{\mathbf{u}}
\newcommand{\qq}{\mathbf{q}}

\newcommand{\mJ}{\mathcal{\partial}}

\newcommand{\tr}{\text{tr}}
\newcommand{\QQ}{\mathbf{q}}

\providecommand{\divv}[1]{\nabla\cdot#1}
\providecommand{\norm}[1]{\lVert#1\rVert}

\providecommand{\keywords}[1]{\textbf{\textit{Index terms---}} #1}

\newtheorem{lemma}{Lemma}
\newtheorem{theorem}{Theorem}
\begin{document}

\title{Iterative solvers for Biot model under small and large deformation} 



\author{  Manuel Borregales\footnote{Department of Mathematics, University of Bergen, PO Box 7800, Bergen, Norway; Emails: \{Manuel.Borregales,  Jan.Nordbotten, Florin.Radu\}\texttt{@uib.no}} \and Kundan Kumar\footnote{Department of Mathematics and Computer Science, Karlstad University, 651 88 Karlstad, Sweden; Email: Kundan.Kumar\texttt{@kau.se}} \and Jan Martin Nordbotten$^{*}$ \and Florin Adrian Radu$^{*}$ 
}

\date{30/05/2019}

\maketitle

\begin{abstract}
We consider $L$-scheme and Newton based solvers for Biot model under small or large deformation. The mechanical deformation follows the Saint Venant-Kirchoff constitutive law. Further, the fluid compressibility is assumed to be nonlinear. A Lagrangian frame of reference is used to keep track of the deformation. We perform an implicit discretization in time (backward Euler) and propose two linearization schemes for solving the nonlinear problems appearing within each time step: Newton’s method and $L$-scheme. The linearizations are used monolithically or in combination with a splitting algorithm. The resulting schemes can be applied for any spatial discretization. The convergences of all schemes are shown analytically for cases under small deformation. Illustrative numerical examples are presented to confirm the applicability of the schemes, in particular, for large deformation.
\end{abstract}

\section{Introduction}
\label{intro}

The coupling of flow and mechanics in a porous medium, typically referred to as poromechanics, plays a crucial role in many socially relevant applications. These include geothermal energy extraction, energy storage in the subsurface, $CO_2$ sequestration, and understanding of biological tissues. The increased role played by computing in the development and optimisation of (industrial) technologies for these applications implies the need for improved mathematical models in poromechanics and robust numerical solvers for them. 

The most common mathematical model for coupled flow and mechanics in porous media is the linear, quasi-stationary Biot model \cite{Biotsttlement,Biot,biot1954,SHOWALTER_2016}. The model consists of two  coupled partial differential equations, representing balance of forces for the mechanics and conservation of mass and momentum for (single-phase) flow in porous media.

In terms of modelling, Biot's model has been extended to unsaturated flow \cite{Jakub_2017,List2016}, multiphase flow \cite{hong,Kraus_2019,mardal,Fullycouple4,RaduNPK15}, thermo-poro-elasticity \cite{Coussy1989}, and reactive transport in porous media \cite{Radu_2013_transport,Radu_2010}, where nonlinearities arise in the flow model, specifically  in the diffusion term, the  time derivative term  and/or in Biot's coupling term. The mechanics model can also be extended to the elasto-plastic  \cite{Explicit1,Ganis_2017}, the fracture propagation \cite{mikelic2013phase} and the hyperelasticity \cite{coussy1995mechanics,coussy},  where the nonlinearities appear in the constitutive law of the material, in the compatibility condition and/or the conservation of momentum equation. 
  Furthermore,  elastodynamics or non-stationary Biot, i.e. Biot-Allard model \cite{Mikelic3}, includes a convolution in the coupling term of both mechanics and flow equations. In this paper, we are going to explore a general case that allows large deformations. The mechanical deformation follows the Saint Venant-Kirchoff constitutive law and the fluid compressibility in the fluid equation is assumed to be nonlinear. This model formulation is needed to later consider extensions of Biot's model to plasticity, more general hyperelastic materials, and elastodynamics.  
 
Finding closed-form solutions for coupled problems is very difficult,  and commonly based on various simplifications. We, therefore, resort to numerical approximations. In general, there are two approaches to solve such problems, the
fully coupled   and the weakly coupled scheme. In general the fully coupled schemes for fluid potential and mechanical deformation are stable, have excellent convergence properties, and ensure that the numerical solution is consistent with the underlying continuous differential equations \cite{Fullycoupled3,Fullycoupled1}. 
Despite obvious advantages, the monolithic solver for the fully coupled problem are more difficult to implement, and have difficulties solving the  resulting linear system, particularly in the context of existing legacy codes for separate physics.
In the weakly coupled approach, while marching in time, we time-lag the flow problem (or the mechanics), thereby fully decoupling the two problems. 
Due to the complexities associated with the fully coupled scheme, the industry standard remains to use weakly coupled or iteratively coupled approaches \cite{Chin,Petersen2012,Settari_2001,Explicit2}. An iteratively coupled approach takes somewhat of a middle path; at each time step, it decouples the flow and mechanics, but iterates so that the convergence is achieved.  Weakly coupled schemes, wherein there are not iterations within time step, have in particular been questioned in previous works \cite{Ferronato_2016,FlorianFCP,Petersen2012,Prevost2013}; they have been shown to lack robustness and even convergence, if not properly designed. In order to ensure the robustness and accuracy of the resulting computations, it is therefore essential to understand the efficiency, stability, and convergence of iterative coupling schemes, in particular in the presence of nonlinearities. 

In this work, we present monolithic and splitting approaches for solving this nonlinear system, {that is, nonlinear compressibility and the Saint Venant-Kirchoff constitutive law for stress-strain}. {Moreover, we rigorously study the convergence of our schemes, including the Newton based ones, under the assumption of small deformations}.  As for splitting approach, we use the undrained split method, see \cite{Kim20112094,Mikelic2}.  We use linear conformal Galerkin elements for the discretization of the mechanics equation and mixed finite elements for the flow equation \cite{Berger2017,GayX2,Jha2007,wheeler,yi2016}. Precisely, the lowest order Raviart-Thomas elements are used \cite{brezzi2012mixed}. We expect, however, that the solution strategy discussed herein will be applicable to other combinations of spatial discretizations such as those discussed in \cite{Jan_2016,Carmen_2016} and the references therein. Backward Euler is used for the temporal discretization.

To summarise, the new contributions of this paper are
\begin{itemize}
\item[$\bullet$]{We propose Newton and $L$-scheme based monolithic and  splitting schemes for solving the Biot model under small or large deformation.}
\item[$\bullet$]{The convergence  analysis of all schemes is  shown  rigorously under the assumption of small deformations.}
\item[$\bullet$]{We provide a benchmark for the convergence of splitting algorithms for a general nonlinear Biot model that includes large deformations.}
\end{itemize}

{We mention some relevant works in this direction}. For the convergence analysis of the undrained split method applied to the linear Biot model, we refer to \cite{Bause_2018,Uwe3,Jakub_2016,GASPAR2017526,Girault,Mikelic2}. For a discussion on the stabilization/tuning parameter used in the undrained split approach, we refer to \cite{Jakub_2016,Jakub_2018}. A theoretical investigation on the optimal choice for this parameter is performed in \cite{Erlend}. The linearization is based on either Newton's method, or the $L$-scheme \cite{List2016,RaduPopKnabner2004,RaduNPK15} or a combination of them \cite{Jakub_2017,List2016}. For monolithic and splitting schemes based solely on $L$-scheme, we refer to \cite{paper_A}. Multirate time discretizations or higher order space-time Galerkin method has also been proposed for the linear Biot model in \cite{Kundan_2016} and \cite{Uwe3}, respectively. 

The paper is structured as follows. In the next section, we present the mathematical model. In Section \ref{IterativeSchemes}, we propose four iterative schemes. Section \ref{secc:SmalDef} shows the  analysis of iterative schemes under the assumption of small deformations. Numerical results are presented in Section \ref{Numerical_Results} followed by the conclusion in Section \ref{Conclusions}. 

\section{Governing equations}
\label{sec:1}

We consider a fluid flow problem in a poroelastic bounded reference domain  $\Omega \subset \mathbb{R}^d$, $d \in \{ 2, 3 \}$  under large deformation. 
A Lagrangian frame of reference is used to keep track of the invertible transformation $x:=\{x(X,t)=X + \uu(X,t): X\in \Omega \rightarrow x \in  \Omega_t \}$, where $\Omega_t$ is the deformed domain at time $t$ and $\uu$ represents the deformation field. The gradient of the transformation and its determinant are given by $\bF = \nabla\, x(X,t)$ and  $J =  \det (\bF)$. All differentials are with respect to the undeformed coordinates $X$, unless otherwise stated. 

We will now write the conservation of momentum and mass equation in $\Omega$.
 The conservation of momentum represents  the balance between the first Piola-Kirchhoff poroelastic stress ${\bf \Pi}$ in $\Omega$ and the forces acting on $\Omega_t$, and is given by  
 
\begin{equation}
-\nabla \cdot {\bf \Pi}   = \rho_b \ga,
\end{equation}
where $ \rho_b = J \varrho_{b} $ is  the bulk density in $\Omega$, $\varrho_{b}$ is the bulk density in $\Omega_t$  and $\ga$ is gravity.

We exploit the relation ${\bf \Pi} =  \bF {\bf \Sigma}$ since the constitutive laws are developed for the second Piola-Kirchhoff poroelastic stress $\bf \Sigma$. This stress tensor is composed of the effective mechanical stress $ {\bf \Sigma}^{eff}$ and the pore pressure $p$ by the following relation
  $${\bf \Sigma} = {\bf \Sigma}^{eff} - J \bF^{-1} \bF^{\top} p,$$
where $J \bF^{-1} \bF^{\top}$ ensures that pressure $p$ exerts an isotropic stress in $\Omega_t$.  We assume an isotropic poroelastic material with constant shear modulus   $\mu$ and a nonlinear function of the volumetric strain $\cc(\cdot)$ \cite{paper_A,Temam_Miranville}. The effective stress is given by Saint Venant-Kirchhoff constitutive law:
 $\bf \Sigma^{eff} =2\mu {\bf E}+ \cc \left( \tr({\bf E}) \right), $ where the Green strain tensor ${\bf E}$ is defined by\\  ${\bf E} = \frac{1}{2}\left( \nabla \uu +  \nabla^{\top} \uu +  (\nabla \uu)^{\top}   \nabla \uu\right)$.

The conservation of fluid mass is given by 
\begin{equation}
 \dot{{\Gamma}} + \divv  \QQ = S_{f}.
\end{equation}
\noindent We consider a  fluid mass  ${\Gamma}  = J  \rho_f \phi$ of a slightly compressible fluid, where $\phi$ is the  porosity and $\rho_f$ the fluid density and $S_{f} $ the source term in $\Omega$ respectively. 
%
The time derivative of the fluid content $\dot{\Gamma} = \dot{\Gamma}(\uu,p)$ is considered to be a function of the pressure and the pore volume change due to the deformation field. We consider Darcy's law 
\begin{equation}
 \QQ = -\bK(\uu) \left( \nabla p -\rho_f \ga_0 \right),
 \end{equation}
where the flux variable $\qq$ is the first Piola transform of the corresponding flux variable in $\Omega_t$, $\bK = J \bF^{-1} {\bf{k}} \bF^{-\top}$ is the corresponding transformation of the mobility tensor  ${\bf{k}}$  in $\Omega_t$ and  $\Upsilon = \bF^{\top} \ga$. Finally, the general nonlinear Biot model considered in this paper reads as:

Find $\left(\uu,\qq,p\right)$ such that
\begin{align}
-\nabla \cdot {\bf \Pi} \left(\nabla \uu ,  p \right) & = \varrho_{b} \ga ,  & \text{in}\ \Omega \times ]0, T [, \nonumber\\
 \QQ &= -\bK(\uu) \left( \nabla p -\rho_f \Upsilon \right),  & \text{in}\ \Omega \times ]0, T [, \label{darcy}\\
\dot{\Gamma}(\uu,p) + \divv \QQ &= S_{f},  & \text{in}\ \Omega \times ]0, T [. \nonumber
\end{align}
\noindent To complete the model we consider  Dirichlet boundary conditions (BC) and initial conditions given by $(\uu_0,p_0)$ such that $\Gamma(\uu_0,p_0) = \Gamma_0$ and ${\bf \Pi}(\uu_0,p_0) = {\bf \Pi}_0$ at time $t = 0$. The functions $\Gamma_0$ and ${\bf \Pi}_0$  are supposed to be given (and to be sufficiently regular). 

In practice, the initial data $\uu_0$ and $p_0$ are not independent and can be obtained by solving the flow equation for $p_0$ and then solving the mechanics equation for getting $\uu_0$.

\section{Iterative schemes}
\label{IterativeSchemes}
In this section, we present several monolithic and splitting iterative schemes for solving Eqs. \eqref{darcy}. First, we propose the Newton method which is well known for having quadratic convergence. Secondly, we combine the Newton method with a stabilized splitting scheme based on the undrained split method. Finally, for the third and fourth schemes, we propose monolithic and splitting $L$-schemes. The iterative schemes will be written using an incremental formulation. In this regard, we introduce naturally defined residuals for the nonlinear Eqs. \eqref{darcy}.

\begin{align}
\begin{array} {ll}
\mathcal{F}_{\rm{mech}} (\uu,p)&= -\nabla \cdot {\bf \Pi} \left( \nabla \uu ,  p \right) - \rho_{b} \ga,  \\  [1.3ex]
\mathcal{F}_{\rm{darcy}}(\uu,p) &= \QQ  + \bK(\uu) \left( \nabla p -\rho_f \Upsilon \right), \\  [1.3ex]
\mathcal{F}_{\rm{mass} }(\uu,p)&= \dot{\Gamma}(\uu,p) + \divv \QQ - S_{f}.
\end{array}
\end{align}

\noindent We will denote by $\delta (\cdot)^{\ii} = (\cdot)^{\ii}-(\cdot)^{\iim} $ the incremental operator, $\ii$ the incremental counter, $\partial_{(\cdot)} $ the partial derivative operator respect to $(\cdot)$.

\subsection{\bf A monolithic Newton solver}

The Newton method is usually the first choice of the linearization methods due to its quadratic convergence. However, the convergence is local and it requires relatively small time steps to ensure the quadratic convergence \cite{Paper_Newton}.  The method starts by using  initial solution $\left( \uu^0,  \qq^0,  p^0 \right) $, solves for $\left( \delta \uu^{\ii},  \delta \qq^{\ii}, \delta p^{\ii} \right)$ satisfying
  
{\fontsize{8}{11}\selectfont  
\begin{align}
\label{Monolithic_N}
\begin{array} {ll}
-\nabla \cdot  \left( \mJ_{\uu} {\bf \Pi} \left( \nabla \uu^{\iim} ,  p^{\iim} \right)  \nabla \delta\uu^{\ii} -  \mJ_{p} {\bf \Pi} \left( \nabla \uu^{\iim} ,  p^{\iim} \right) \delta p^{\ii}    \right) &=  -\mathcal{F}_{mech} (\uu^{\iim},p^{\iim}),  \\  [1.3ex]
\delta \QQ^{\ii}+ \bK(\uu^{\iim})\nabla \delta p^{\iim} + \mJ_{\uu} \bK(\uu^{\iim})\nabla p^{\iim} \delta \uu^{\ii}& = -\mathcal{F}_{darcy} (\uu^{\iim},p^{\iim}),\\  [1.3ex]
\mJ_p \dot{\Gamma}(\uu^{\iim},p^{\iim}) \delta p^{\ii}+ \mJ_{\uu} \dot{\Gamma}(\uu^{\iim},p^{\iim}) \delta {\uu}^{\ii} + \divv \delta\QQ^{\ii} &=-\mathcal{F}_{mass} (\uu^{\iim},p^{\iim}), 
\end{array}
\end{align}
}
and finally updates the variables $$\left(  \uu^{\ii},   \qq^{\ii},  p^{\ii} \right)  = \left(  \uu^{\iim},   \qq^{\iim},  p^{\iim} \right) + \left( \delta \uu^{\ii},  \delta \qq^{\ii}, \delta p^{\ii} \right). $$

\subsection{{\bf A splitting Newton solver}}

The splitting Newton method combines a splitting method with the Newton linearization. 
We introduce a stabilization  parameter $L_s \ge 0$ to stabilize the mechanics equation. The precise condition on $L_s$ to ensure convergence  is shown in Theorem \ref{Newton_Spliting_convergence}.
The  method consists on two steps: starting with the initial  condition $ \left(\uu^0,  \qq^0,  p^0 \right) $:

\textbf{Step 1:} solve for $\left( \delta \QQ^{\ii} ,\delta p^{\ii} \right) $

\begin{align}
\label{Splitting_N}
\begin{array} {ll}
\delta \QQ^{\ii} + \bK(\uu^{\iim})\nabla \delta p^{\iim}& = -\mathcal{F}_{darcy} (\uu^{\iim},p^{\iim}),   \\  [1.3ex]
\mJ_p \dot{\Gamma}(\uu^{\iim},p^{\iim}) \delta p^{\ii} + \divv \delta\QQ^{\ii} &=-\mathcal{F}_{mass} (\uu^{\iim},p^{\iim}) ,
\end{array}
\end{align}

and update the variables $$\left(  \qq^{\ii},  p^{\ii} \right)  = \left(   \qq^{\iim},  p^{\iim} \right) + \left( \delta \qq^{\ii}, \delta p^{\ii} \right). $$

\textbf{Step 2:} solve for $\delta \uu^{\ii}$ satisfying

\begin{align}
-\nabla \cdot \left(  \mJ_{\uu} {\bf \Pi} \left(\nabla \uu^{\iim} ,  p^{\ii} \right)  \nabla \delta \uu^{\ii} - L_s  (\nabla \cdot  \delta\uu^{\ii}) \ \bI  \right)  &=  -\mathcal{F}_{mech} (\uu^{\iim},p^{\ii}) , \label{Splitting_N_mech} 
\end{align}

and update the variable $$  \uu^{\ii}=  \uu^{\iim} + \delta \uu^{\ii}.$$
The stability of the scheme is controlled by $L_s $ as it is shown in \cite{Paper_Newton}.

\subsection{{\bf A monolithic $L$-scheme}}

The $L$-scheme can be interpreted as either a stabilized Picard method or a quasi-Newton method. This scheme is  robust but only linearly convergent. Moreover, it can be applied to non-smooth but monotonically increasing nonlinearities. For example, for the case of  H\"older continuous (not Lipschitz) nonlinearities we refer to \cite{Jakub_Sorin_2018}. As it is a fixed point scheme, it  can be speeded up by using the Anderson acceleration \cite{anderson,Jakub_2018}. To summarize, the main advantages of the $L$-scheme are:
\begin{itemize}
\item{It does not involve computation of derivatives.}
\item{The arising linear systems are well-conditioned. }
\item{It can be applied to non-smooth nonlinearities.}
\item{It is easy to understand and implement.}
\end{itemize}
\noindent A monolithic $L$-scheme requires three constant tensors ${\bf{L}}_{\uu},\ {\bf{L}}_{p},\ {\bf{L}}_{\qq} \in \mathbb{R}^{d\times d}$  and two positive constants $L_p $ and  $L_{\uu} $ as linearization parameters. A practical choice of the linearization parameters will be discussed in the numerical section. We refer to \cite{paper_A,FlorianFCP} for a discussion regarding the best  choice for the linearization parameters $L_p $ and $L_{\uu}.$

  The method starts with the given initial solution $\left( \uu^0,  \qq^0,  p^0 \right) $ and solve for $\left( \delta \uu^{\ii},  \delta \qq^{\ii}, \delta p^{\ii} \right) $ 


\begin{align}
\begin{array}{ll}
-\nabla \cdot {\bf{L}}_{\uu}  \nabla \delta \uu^{\ii} - \nabla \cdot {\bf{L}}_{p} \delta p^{\ii}     &=  -\mathcal{F}_{mech} (\uu^{\iim},p^{\iim}), \\  [1.3ex]
\delta \QQ^{\ii} + \bK(\uu^{\iim})\nabla \delta p^{\ii} + {\bf{L}}_{\qq}  \delta \uu^{\ii}& = -\mathcal{F}_{darcy} (\uu^{\iim},p^{\iim}),\\  [1.3ex]
L_p \delta p^{\ii}+ L_\uu \delta {\uu}^{\ii} + \divv \delta\QQ^{\ii} &=-\mathcal{F}_{mass} (\uu^{\iim},p^{\iim}),
\end{array}
\end{align}

and then update the variables $$\left(  \uu^{\ii},   \qq^{\ii},  p^{\ii} \right)  = \left(  \uu^{\iim},   \qq^{\iim},  p^{\iim} \right) + \left( \delta \uu^{\ii},  \delta \qq^{\ii}, \delta p^{\ii} \right). $$

\subsection{{\bf A splitting $L$-scheme}}

The splitting scheme requires less linearization terms: two constants  $\mathbf{L}_{\uu}\in \mathbb{R}^{d\times d}$, $L_p\ge 0$ and a positive stabilisation term $L_s$. This makes it suitable for quick implementation since there is no need to calculate any Jacobian.
  The  method is split in two steps, given initial solution $ \left(\uu^0,  \qq^0,  p^0 \right) $:

\textbf{Step 1:} solve for $\left( \delta \QQ^{\ii} ,\delta p^{\ii} \right) $

\begin{align}
\begin{array}{ll}
\delta \QQ^{\ii} + \bK(\uu^{\iim})\nabla \delta p^{\ii}& = -\mathcal{F}_{darcy} (\uu^{\iim},p^{\iim}), \\  [1.3ex]
L_p \delta p^{\ii} + \divv \delta\QQ^{\ii} &=-\mathcal{F}_{mass} (\uu^{\iim},p^{\iim}) ,
\end{array}
\end{align}

update the variables $$\left(  \qq^{\ii},  p^{\ii} \right)  = \left(   \qq^{\iim},  p^{\iim} \right) + \left( \delta \qq^{\ii}, \delta p^{\ii} \right). $$

\textbf{Step 2:} solve for $\delta \uu^{\ii}$
\begin{align}
-\nabla \cdot \left(  {\bf{L}}_{\uu}    \nabla \delta \uu^{\ii} +{L_s}  (\nabla \delta \cdot \uu^{\ii})  \ \bI \right)&=  -\mathcal{F}_{mech} (\uu^{\iim},p^{\ii}),   
\end{align}

and then  update the variables $$  \uu^{\ii}=  \uu^{\iim} + \delta \uu^{\ii}.$$

\section{The Biot model under small deformations}
\label{secc:SmalDef}
 The convergence analysis of the iterative schemes proposed cannot be addressed with standard techniques \cite{paper_A,Jakub_2018,Jakub_2017,List2016,Mikelic2}.  This is due to the nonlinearities being non-monotone. Nevertheless, a rigorous analysis can be performed for the case of small deformations.
Accordingly, we assume the porous medium to be under small deformation and present the convergence of the iterative schemes proposed in the previous section.

Under small deformation, the different between  $\Omega_t$ and $\Omega$ can be neglected. The gradient of the transformation is approximated by $\bF \approx \bI $ and the determinant of the transformation by $J \approx 1$.  Additionally, the Green strain tensor $\bf E$ can be approximated by the infinitesimal strain tensor ${\bf E} \approx \varepsilon = \frac{1}{2} \left( \nabla \uu + (\nabla \uu)^{\top}\right)$.
Then, the poroelastic stress tensor can be expressed by 
\begin{equation}
{\bf \Pi }(\uu,p) =\sigma(\uu,p)= 2\mu \varepsilon(\nabla \uu) +  \mathfrak{c}(\tr (\varepsilon(\nabla \uu))) - \alpha p \bI, \label{small_Stress}
\end{equation}
where $\alpha$ is the Biot constant. The mobility tensor is considered isotropic ${\bf K} (\uu,p)  = k \bI $, but the results of the convergence analysis can be extended  without difficulties to a more general anisotropic case. Additionally, the time derivative of the volumetric deformation is approximated by $\dot{J} \approx \divv \dot{\uu} $. In this regard the fluid mass can be expressed as
\begin{align}
\Gamma (\uu,p) & = \Gamma_0 + c_f \left( \bb(p)-\bb(p_0) \right) + \alpha \divv \left(\uu - \uu_0 \right), \label{small_mass}
\end{align}
 where the relative density $\bb (\cdot)$ is a nonlinear function of the pressure $p$. The variational formulation for the Biot model, under small deformation, reads as follows:

For each $t\in (0,T]$, find ${\mathbf u}(t) \in \left(H_0^1(\Omega)\right)^d$, $\qq \in H^1( \rm{div},\Omega) $ and $p(t) \in  L^2(\Omega) $ such that there holds

\begin{align}
\label{fully1}
\begin{array}{rl}
\left(\varepsilon ({\mathbf u}), \varepsilon ({\mathbf v})\right) + \left(\cc(\nabla \cdot {\mathbf u})-\alpha p,
\nabla \cdot {\mathbf v}\right)  &= (\rho_b \ga, {\mathbf v}), \ \ \forall {\mathbf v} \in \left(H(\Omega)\right)^d,  \\  [1.3ex]
\left(\bK^{-1} \qq , {\bf z} \right) - \left( p ,\divv {\bf z} \right) &= \left(\rho_f \ga ,{\bf z} \right), \ \ \forall {\bf z} \in H^1(\text{div},\Omega), \\  [1.3ex]
\left( \dot{\bb}(p)) + \alpha   \nabla \cdot \dot{{\mathbf u}} , w \right)  
+ \tau \left( \divv \qq, w \right) &= \tau \left( S_f,w\right), \ \forall w \in L^2(\Omega), 
\end{array}
\end{align}

with the initial condition 
 \begin{equation}
 \left( \bb(p_0)) + \alpha   \nabla \cdot {\mathbf u}_0 , w \right) = 0, \  \forall w \in L^2(\Omega).
 \end{equation}

\noindent In the above, we have used the standard notations. We denote by $L^2(\Omega)$ the space of square integrable functions and by $H^1(\Omega)$ the Sobolev space $H^1(\Omega) = \{v \in L^2(\Omega)\,;\, \nabla\,v \in L^2(\Omega)^d\}.$
Furthermore,  $H^1_0(\Omega) $ will be the space of functions in $H^1(\Omega)$ vanishing on $\partial \Omega$ and $H({\rm div};\Omega)$  the space of vector valued function having all the components and the divergence in $L^2(\Omega)$.
As usual we denote by $(\cdot,\cdot)$ the inner product in $L^2(\Omega)$, and by $||\cdot||$ its associated norm.

Next, we make structural assumptions on the nonlinearities:
\begin{itemize}
\item[(A1)] $\cc,\, \bb:{\mathbb R} \rightarrow {\mathbb R}$ differentiable with $\cc'$ and $\bb'$ Lipschitz continuous.
\item[(A2)] There exists a constant $\alpha_{\cc}$ such that $\cc'(\xi)>~\alpha_{\cc}$, $\forall \, \xi \in {\mathbb R}$.
\item[(A3)] There exists a constant $\alpha_{\bb}$ such that $\bb'(\xi)>\alpha_{\bb}$, $\forall \, \xi \in {\mathbb R}$.
\item[(A4)] There exists  constant $k_{m}>0$ and $k_{M}$ such that  
$ k_{m}  \leq k(\vec{\xi})\leq  k_{M},$ $\forall {\vec{ \xi}} \in \Omega$.
\end{itemize}

For the discretization of problem \eqref{fully1} we use conformal Galerkin finite elements for the displacement variable and mixed finite elements for the flow  \cite{GayX2,wheeler}. More precisely, we use linear elements for the displacement and lowest order Raviart-Thomas elements \cite{brezzi2012mixed} for the flow. Backward Euler is used for the temporal discretization.

Let  $\Omega = \cup_{K \in \mathcal{T}_h} K$ be a regular decomposition of $\Omega$ into $d$-simplices. We denote by $h$ the mesh size. The discrete spaces are given by
\begin{align}
\hspace{1cm}{\mathbf V}_h &:= \{\mathbf{v}_h \in {H^1(\Omega)}^d\,;\,  {{\mathbf{v}}_h}_{|K} \in {{\mathbb{P}_1^d}} \,, \,   {\forall} K \in  \mathcal{T}_h\}, \nonumber\\
\hspace{1cm}W_h &:= \{w_h \in L^2(\Omega) \,;\, {w_h}_{|K} \in {{\mathbb{P}}_0} \, ,\,  {\forall} K \in  \mathcal{T}_h \}, \nonumber\\
 \hspace{1cm} {\mathbf{Z}}_h &:= \{{\vec{z}}_h \in {H(\mathrm{div}; \Omega)}\,;\,  {\vec{z}}_{h|K} (\vec{x})=\vec{a}+b \vec{x}, \, \vec{a} \in  \mathbb{R}^d,\, b \in \mathbb{R} ,  \, {\forall} K \in  \mathcal{T}_h\},\nonumber
\end{align}
where ${\mathbb{P}}_0, {\mathbb{P}}_1$ denote the spaces of constant functions and of linear polynomials, respectively.  For $N \in \mathbb{N}$, we discretize the time interval uniformly and define the time step $\tau = \frac{T}{N}$ and $t_n = n\tau$. We use the index $n$ for the primary variable $\uu^{n}$, $\qq^{n}$ and $p^{n}$  at corresponding time step $t_n$.
In this way, the fully discrete weak
problem reads:
 
For $n\geq 1$ and given $\left( \uu_h^{n-1}, \qq_h^{n-1},p_h^{n-1} \right)$ find $\left( \uu_h^n, \qq_h^n,p_h^n \right) \in \left( {\mathbf V}_h,{\mathbf Z}_h,W_h \right) $, such that

\begin{align}
\label{fully1_discret}
\begin{array}{rl}
\left(\varepsilon ({\mathbf u}_h^{n}), \varepsilon ({\mathbf v}_h)\right) + \left(\cc(\nabla \cdot {\mathbf u}_h^{n}),\nabla \cdot {\mathbf v}_h\right)  - \alpha \left(  p_h^{n}, \nabla\cdot {\mathbf v}_h\right) &= (\rho_b \ga, {\mathbf v}_h),   \\ [1.3ex]
\left(\bK^{-1} \qq_h^{n} , {\bf z}_h \right) - \left( p_h^{n} ,\divv{\bf z}_h \right) &= \left(\rho_f \ga ,{\bf z}_h \right),  \\ [1.3ex]
\left( \bb(p_h^{n}) - \bb(p_h^{n-1}), w_h \right) + \alpha \left(    \nabla \cdot ({\mathbf u}_h^{n} -  {\mathbf u}_h^{n-1}), w_h \right)    \\  [1.3ex]
+ \tau \left( \divv \qq_h^{n},\nabla w_h \right) &= \tau (S_f,w_h),
\end{array}
\end{align}

$\text{for all } \left( {\bf v}_h, {\bf z}_h, w_h \right) \in \left( {\mathbf V}_h,{\mathbf Z}_h,W_h \right) $.

Following the notation previously introduced, we denote by $n$ the time level, whereas $i$ will refer to the iteration number of the Newton method. We further denote the approximate solution of the linearized problem~\eqref{fully1_discret} by $({\bu}_h^{n,i}, \qq_h^{n,i}, p_h^{n,i})$. At this stage we can introduce the notations
\begin{eqnarray*}
{\be}_{\bu}^{n,i} &=& {\bu}_h^{n,i}-{\bu}_h^{n},  \\[1ex]
{\be}_{\qq}^{n,i} &=& {\qq}_h^{n,i}-{\qq}_h^{n},  \\[1ex]
e_p^{n,i} &=& p_h^{n,i} - p_h^{n}.
\end{eqnarray*}
These will be used subsequently in the convergence analysis of the  monolithic Newton method and the alternate version. For the monolithic and splitting $L$-scheme  the convergence analysis can be found in \cite{paper_A}.

\subsection{\bf Convergence analysis of the monolithic Newton method}\label{subsec:monolithic}

In this section, we analyse the monolithic Newton method introduced in Section \ref{IterativeSchemes} used for solving the simplified nonlinear Biot model given in~\eqref{fully1_discret}. As we have previously stated, we perform the analysis for the case of small deformation. Here we present a variational formulation of the scheme and demonstrate its quadratic convergence in a rigorous manner. 
 The Newton scheme reads as follows:\\ [1.3ex]
\noindent For $i=1,2,\ldots$ solve\\ [-2ex]

{\fontsize{9}{11}\selectfont 
\begin{align}
\label{small_fully1_Newton}
\begin{array}{rl}
\left(\varepsilon ({\mathbf u}_h^{n,\ii}), \varepsilon ({\mathbf v}_h)\right) + \left(\cc(\nabla \cdot {\mathbf u}_h^{n,\iim})+\cc'(\nabla \cdot {\mathbf u}_h^{n,\iim})\nabla \cdot  \delta {\mathbf u}_h^{n,\ii},
\nabla \cdot {\mathbf v}_h\right)  \\ [1.3ex]
 - \left( \alpha p_h^{n,\ii}, \nabla\cdot {\mathbf v}_h\right) = (\rho_b \ga, {\mathbf v}_h),  \\ [1.3ex]
\left(\bK^{-1} \qq_h^{n,\ii} , {\bf z}_h \right) - \left( p_h^{n,\ii} ,\divv {\bf z}_h \right) = \left(\rho_f \ga ,{\bf z}_h \right), \\ [1.3ex]
\left( \bb(p_h^{n,\iim})+\bb'(p_h^{n,\iim})\delta p_h^{n,\ii} - \bb(p_h^{n-1}), w_h \right) + \left( \alpha   \nabla \cdot ({\mathbf u}_h^{n,\ii} -  {\mathbf u}_h^{n-1}), w_h \right)    \\ [1.3ex]
+ \tau \left( \divv \qq_h^{n,\ii},\nabla w_h \right) = \tau (S_f,w_h),  
\end{array}
\end{align} 
}
 $\forall \left( {\mathbf v}_h,{\mathbf z}_h , w_h\right)\in \left( {\mathbf V}_h,{\mathbf Z}_h,W_h\right)$,
where the initial approximation $({\mathbf u}_h^{n,0} \qq_h^{n,0},p_h^{n,0})$ is taken as the solution at the previous time step, that is $({\mathbf u}_h^{n-1}, \qq_h^{n-1},p_h^{n-1})$.\\

In order to prove the convergence of the considered Newton method, the following lemmas will be used.

\begin{lemma} \label{succeion_convergence}
Let $\left\lbrace x_n \right\rbrace_{n\geq0}$  be a sequence of real positive number satisfying 
\begin{equation}
x_n \leq a x_{n-1}^{2} + b x_{n-1} \ \forall n \geq 1,
\end{equation}
where $a,b \ge 0$. Assuming that
 $$ a x_0^2 + b \leq 1$$ 
 holds, then the sequence $\left\lbrace x_n \right\rbrace_{n\geq0}$ converges to zero.
\end{lemma}

\proof
The result can be shown by induction, see page 52 in \cite{Florin_thesis} for more details.
\endproof

\begin{lemma}\label{lemma1}
If $f:{\mathbb R} \rightarrow {\mathbb R}$ is differentiable and $f'$ is Lipschitz continuous, then there holds
$$|f(x)-f(y)+ f'(y)(y-x)|\leq \displaystyle \frac{L_{f'}}{2}|y-x|^2, \quad \forall \, x,y\in{\mathbb R}.$$
\end{lemma}
\proof
See page $350$ in~\cite{Knabner_2003}, for example.
\endproof
\noindent Next, the following result provides the quadratic convergence of the Newton method~\eqref{small_fully1_Newton} for $\tau$ sufficiently small.

\begin{theorem}
Assuming (A1)-(A4), the Newton method in~\eqref{small_fully1_Newton} converges quadratically if $\tau = O (h^d)$. 


\label{Theorem1}
\end{theorem}

\proof
 By subtracting equations~\eqref{fully1_discret} from~\eqref{small_fully1_Newton},  taking as test functions ${\mathbf e}_{\mathbf u}^{n,\ii}$, ${\be}_{\qq}^{n,\ii}$ and $e_p^{n,\ii}$ and rearranging some terms to the right hand side we obtain,

{
\begin{align}
\left(\varepsilon ({\mathbf e}_{\uu}^{n,\ii}), \varepsilon ({\mathbf e}_{\uu}^{n,\ii})\right) + \left( \cc'(\nabla \cdot {\mathbf u}_h^{n,\iim})\nabla \cdot  \be_{\uu}^{n,\ii},
\nabla \cdot {\mathbf e}_{\uu}^{n,\ii} \right) - \alpha \left(  e_p^{n,\ii}, \nabla\cdot {\mathbf e}_{\uu}^{n,\ii}\right)
\nonumber \\
 = \left(\cc(\nabla \cdot {\mathbf u}_h^{n})-\cc(\nabla \cdot {\mathbf u}_h^{n,\iim})+\cc'(\nabla \cdot {\mathbf u}_h^{n,\iim})\nabla \cdot   \be_{\uu}^{n,\iim},
\nabla \cdot {\mathbf e}_{\uu}^{n,\ii} \right), \label{fully1_error} 
 \\ 
\left(\bK^{-1}{\bf e}_{\qq}^{n,\ii} , {\bf e}_{\qq}^{n,\ii} \right) - \left( e_p^{n,\ii} , \divv \be_{\qq}^{n,\ii}\right) = 0, \label{fully3_error} 
\\ 
\left(\bb'(p_h^{n,\iim}) \left( p_h^{n,\ii}-p_h^{n}  \right), e_p^{n,\ii} \right) + \alpha \left(    \nabla \cdot \be_{\uu}^{n,\ii} , e_p^{n,\ii} \right)  
+ \tau \left( \divv \be_{\qq}^{n,\ii},e_p^{n,\ii} \right)
\nonumber  \\
 = \left( \bb(p_h^{n,\iim})-\bb(p_h^{n-1})+\bb'(p_h^{n,\iim}) \left( p_h^{n,\iim}-p_h^{n} \right), e_p^{n,\ii} \right),  \label{fully2_error} 
\end{align}
}


where we have rewritten,


\begin{align}
\begin{array}{ll}
 \cc'(\nabla \cdot {\mathbf u}_h^{n,\iim}) \nabla \cdot {\delta \uu}_h^{n,\ii} &= \cc'(\nabla \cdot {\mathbf u}_h^{n,\iim}) \nabla \cdot ({ \mathbf u}_h^{n,\ii}-{\mathbf u}_h^{n,\iim}) \nonumber\\
&=  \cc'(\nabla\cdot {\mathbf u}_h^{n,\iim})(\nabla\cdot {\mathbf u}_h^{n,\ii}-\nabla \cdot {\mathbf u}_h^{n})
\\& -  \cc'(\nabla\cdot {\mathbf u}_h^{n,\iim})(\nabla \cdot {\mathbf u}_h^{n,\iim}-\nabla \cdot {\mathbf u}_h^{n}) \nonumber\\
&=   \cc'(\nabla\cdot {\mathbf u}_h^{n,\iim}) \left( \nabla \cdot {\mathbf e}_{\mathbf u}^{n,\ii} -\nabla \cdot {\mathbf e}_{\mathbf u}^{n,\iim}\right)
,  \label{rewriting_c}
\end{array}
\end{align}
We obtain  an  analogous expression for the term with $\bb'(\cdot)$. From (A1), $\cc(\cdot)$ is differentiable with $\cc'(\cdot)$ Lipschitz continuous, then from Lemma~\ref{lemma1} we have,
\begin{equation}\label{lemma_for_g}
|\cc(x)-\cc(y)+  \cc'(y)(y-x)|\leq \displaystyle \frac{L_{\cc'}}{2}|x-y|^2, \quad \forall \, x,y\in{\mathbb R},
\end{equation}
where $L_{\cc'}$ represents the Lipshitz constant of $\cc'(\cdot)$.
Then, by using Young's inequality $(a,b)\leq \displaystyle \frac{||a||^2}{2\gamma}+\frac{\gamma ||b||^2}{2}$, for $\gamma\ge 0$, and by choosing $x = \nabla\cdot {\mathbf u}_h^{n}$ and $y = \nabla\cdot {\mathbf u}_h^{n,\iim}$ in~\eqref{lemma_for_g}, from~\eqref{fully1_error} we obtain the following bound, for any $\gamma \ge 0$

\begin{align}\label{first_bound_g}
\begin{array}{rl}
||\varepsilon({\mathbf e}_{\mathbf u}^{n,\ii})||^2+\left(\cc'(\nabla\cdot {\mathbf u}_h^{n,\iim})
\nabla \cdot {\mathbf e}_{\mathbf u}^{n,\ii},\nabla \cdot {\mathbf e}_{\mathbf u}^{n,\ii}\right)
-\alpha \left( e_p^{n,\ii},\nabla \cdot {\mathbf e}_{\mathbf u}^{n,\ii}\right) 
\\   [1.3ex]
\leq \displaystyle \frac{L_{\cc'}^2}{8\gamma}||\nabla \cdot {\mathbf e}_{\mathbf u}^{n,\iim}||_{L^4(\Omega)}^4  
+ \frac{\gamma}{2}||\nabla \cdot {\mathbf e}_{\mathbf u}^{n,\ii}||^2. 
\end{array}
\end{align}

\noindent Next, by using the inverse inequality for discrete spaces $||\cdot||_{L^4(\Omega)}\leq C h^{-d/4}||\cdot||$  \cite{Brenner_1991}, (pg. 111)
 the latter reads,
\begin{align}
\label{first_bound_g_v2}
 \begin{array}{rl}
||\varepsilon({\mathbf e}_{\mathbf u}^{n,\ii})||^2+\left(\cc'(\nabla\cdot {\mathbf u}_h^{n,\iim})
\nabla \cdot {\mathbf e}_{\mathbf u}^{n,\ii},\nabla \cdot {\mathbf e}_{\mathbf u}^{n,\ii}\right)
-\alpha \left( e_p^{n,\ii},\nabla \cdot {\mathbf e}_{\mathbf u}^{n,\ii}\right)
 \\  [1.3ex]
\leq C_1 \displaystyle h^{-d}\frac{L_{\cc'}^2}{8\gamma}||\nabla \cdot {\mathbf e}_{\mathbf u}^{n,\iim}||^4 + \frac{\gamma}{2}||\nabla \cdot {\mathbf e}_{\mathbf u}^{n,\ii}||^2. 
\end{array}
\end{align}
\noindent Finally, by using (A2) and choosing $\gamma=\alpha_{\cc}$, we obtain the following inequality,
\begin{equation}\label{final_bound_g}
||\varepsilon({\mathbf e}_{\mathbf u}^{n,\ii})||^2+\frac{\alpha_{\cc}}{2}||\nabla \cdot {\mathbf e}_{\mathbf u}^{n,\ii}||^2
-\alpha \left( e_p^{n,\ii},\nabla \cdot {\mathbf e}_{\mathbf u}^{n,\ii}\right) \leq C_1 \displaystyle  h^{-d}\frac{L_{\cc'}^2}{8\alpha_{\cc}}||\nabla \cdot {\mathbf e}_{\mathbf u}^{n,\iim}||^4.
\end{equation}

\noindent In a similar way, we obtain the following expression from~\eqref{fully2_error},
\begin{equation}\label{final_bound_b}
\tau \left( \divv \be_{\qq}^{n,\ii},e_p^{n,\ii} \right)
+\frac{\alpha_{\bb}}{2}||e_p^{n,\ii}||^2
+\alpha \left(\nabla \cdot {\mathbf e}_{\mathbf u}^{n,\ii}, e_p^{n,\ii}\right) \leq C_2 \displaystyle h^{-d}\frac{L_{\bb'}^2}{8\alpha_{\bb}}||e_p^{n,\iim}||^4.
\end{equation}
\noindent  Adding ~\eqref{final_bound_g},~\eqref{final_bound_b}, and \eqref{fully3_error} multiplied by $\tau$  yields,
{ \fontsize{9}{11}\selectfont
\begin{align}
\label{final_bound}
\begin{array} {ll}
\frac{\alpha_{\cc}}{2}||\nabla \cdot {\mathbf e}_{\mathbf u}^{n,\ii}||^2
+\frac{\alpha_{\bb}}{2}||e_p^{n,\ii}||^2 +\left(\bK^{-1}{\bf e}_{\qq}^{n,\ii} , {\bf e}_{\qq}^{n,\ii} \right) 
&\leq C_1 \displaystyle h^{-d}\frac{L_{\cc'}^2}{8\alpha_{\cc}}||\nabla \cdot {\mathbf e}_{\mathbf u}^{n,\iim}||^4 
 \\  [1.3ex]
 & + C_2 h^{-d}\frac{L_{\bb'}^2}{8\alpha_{\bb}}||e_p^{n,\iim}||^4.
\end{array}
\end{align}
}

\noindent  By defining  $ \alpha_{\cc,\bb}= \min\left( \alpha_{\cc}, \alpha_{\bb},  \frac{\tau}{k_M} \right) $ and $C_{\cc,\bb} = \max \left( \frac{C_1 L_{\cc'}^2}{\alpha_{\cc}} , \frac{C_2 L_{\bb'}^2}{\alpha_{\bb}} \right)$ we can rewrite \eqref{final_bound} as
\begin{equation}\label{final_bound3}
\norm{\nabla \cdot {\mathbf e}_{\mathbf u}^{n,\ii}}^2
+\norm{e_p^{n,\ii}}^2  +\norm{\be_{\qq}^{n,\ii}}^2
\leq \frac{C_{\cc,\bb} h^{-d}}{\alpha_{\cc,\bb}} \left(\norm{\nabla \cdot {\mathbf e}_{\mathbf u}^{n,\iim}}^4
+ \norm{e_p^{n,\iim}}^4 \right).
\end{equation}
Using $\norm{\nabla \cdot {\mathbf e}_{\mathbf u}^{n,0}}\leq C\tau$, $\norm{ {e}_{p}^{n,0}}\leq C  \tau$ (which can be proven) and Lemma \ref{succeion_convergence}, the  quadratic convergence of Newton's method is ensured if 
$$\frac{C_{\cc,\bb} h^{-d}}{\alpha_{\cc,\bb}} \tau^2 \leq 1$$
 which holds true for  $\tau = O (h^{\frac{d}{2}})$.
\endproof

\subsection{\bf Convergence analysis of the alternate splitting Newton scheme}\label{subsec:alternate}
In this section we present the splitting Newton scheme for solving the nonlinear Biot model given in~\eqref{fully1_discret}.  We present the solver in a variational form and demonstrate its linear convergence.

\noindent Let $i \ge 1$, $L_s\ge 0$ and $(\bu_h^{n,\iim},\qq_h^{n,\iim}, p_h^{n,\iim}) \in \left( {\mathbf V}_h,{\mathbf Z}_h,W_h\right)$ be given.

\textbf{Step 1:} find $(\qq_h^{n,\ii}, p_h^{n,\ii}) \in \left({\mathbf Z}_h,W_h\right)$ such that

{\fontsize{9}{11}\selectfont  
\begin{align}
\label{small_splitting2_Newton}
\begin{array}{rl}
\left(\bK^{-1} \qq_h^{n,\ii} , {\bf z}_h \right) - \left( p_h^{n,\ii} ,\nabla \cdot {\bf z}_h \right) &= \left(\rho_f \ga ,{\bf z}_h \right), \\ [1.3ex]
\left( \bb(p_h^{n,\iim})+\bb'(p_h^{n,\iim})\delta p_h^{n,\ii} - \bb(p_h^{n-1}), w_h \right)  + \tau \left( \divv \qq_h^{n,\ii},\nabla w_h \right)    \\  [1.3ex]
+ \alpha \left(    \nabla \cdot ({\mathbf u}_h^{n,\iim} -  {\mathbf u}_h^{n-1}), w_h \right) &= \tau (S_f,w_h), 
\end{array}
\end{align}
}
 $\forall \left( {\mathbf z}_h , w_h\right)\in \left( {\mathbf Z}_h,W_h\right)$

\textbf{Step 2:} find $\bu_h^{n,\ii} \in  {\mathbf V}_h$ such that

\begin{align}
\label{small_splitting1_Newton} 
\begin{array}{rl}
\left(\varepsilon ({\mathbf u}_h^{n,\ii}), \varepsilon ({\mathbf v}_h)\right) + \left(\cc(\nabla \cdot {\mathbf u}_h^{n,\iim})+\cc'(\nabla \cdot {\mathbf u}_h^{n,\iim})\nabla \cdot  \delta {\mathbf u}_h^{n,\ii},
\nabla \cdot {\mathbf v}_h\right)  
\\  [1.3ex]
+ \left(L_s\divv \delta {\mathbf u}_h^{n,\ii},
\nabla \cdot {\mathbf v}_h\right)
 - \alpha \left(  p_h^{n,\ii}, \nabla\cdot {\mathbf v}_h\right) = (\rho_b \ga, {\mathbf v}_h), 
\end{array}
\end{align}

 $\forall  {\mathbf v}_h \in  {\mathbf V}_h$.

\begin{theorem}
\label{Newton_Spliting_convergence}
Assuming (A1)-(A4) and $L_s\geq \frac{\alpha^2}{\alpha_{\bb}}$, the alternate Newton splitting method in~\eqref{small_splitting2_Newton}-\eqref{small_splitting1_Newton} converges linearly if $\tau$ is small enough.
\end{theorem}
\proof The proof is similar to that of Theorem \ref{Theorem1}. Nevertheless, for the sake of completion we give it in Appendix \ref{Appendix_A}.
\endproof

\section{Numerical examples} 
\label{Numerical_Results}

In this section, we present numerical experiments that illustrate  the performance of the proposed iterative schemes. We study two test problems: a 2D academic problem with a manufactured analytical solution, and a 3D large deformation case on a unit cube. All numerical experiments were implemented using the open-source finite element library Deal II \cite{DealII}. For all numerical experiments, a Backward Euler scheme has been used  for the time  discretization. We consider continuous linear Galerkin FE for $\uu$, lowest order of Raviart-Thomas FE and discontinuous Galerkin FE for  $\qq$  and $p$. However, we would like to mention that any stable discretization can be considered instead. For all cases, as stopping criterion for the schemes, we use 

$$ \norm{p^{\ii}-p^{\iim}} +\norm{{\qq}^{\ii}-\qq^{\iim}} +\norm{{\uu}^{\ii}-\uu^{\iim}} \leq 10^{-8}.$$

\noindent {\bf Test problem 1: an academic example for Biot's model under small deformation}

We solve the nonlinear Biot problem under small deformation in the unit-square  $\Omega=(0,1)^2$ and until final time $T = 1$. This test case was proposed in \cite{paper_A} to study the performance of the monolithic and splitting $L$-scheme. We extend the Newton method and the alternate Newton method described in Section \ref{secc:SmalDef}.

Here, we introduce a manufactured right hand side such that the problem admits the following analytical solution
\begin{align}
p(x,y,t)&=tx(1-x)y(1-y), \ \ \qq(x,y,t)=-k\nabla p, \nonumber        \\
u_1(x,y,t)=u_2(x,y,t)&=tx(1-x)y(1-y), \nonumber
\end{align} 
which has homogeneous boundary values for $p$ and $\uu$.

For infinitesimal deformations and rotations, there is no distinction between the reference and the deformed domains.
In this regard, we solve problem \eqref{fully1_discret} using the iterative schemes proposed in Section \ref{secc:SmalDef}. 
The mesh size and the time step are set as $h = \tau = 0.1$.  For this case, all initial conditions are zero. 
 The linearization parameters $L_p$ and $L_{\uu}$ are equal to the Lipschitz constant $L_{\bb}$ and $L_{\cc}$ corresponding to the nonlinearities $\bb(\cdot)$ and $\cc(\cdot)$  \cite{paper_A}.

\begin{figure}[h!]
\centering
\includegraphics[scale=0.25]{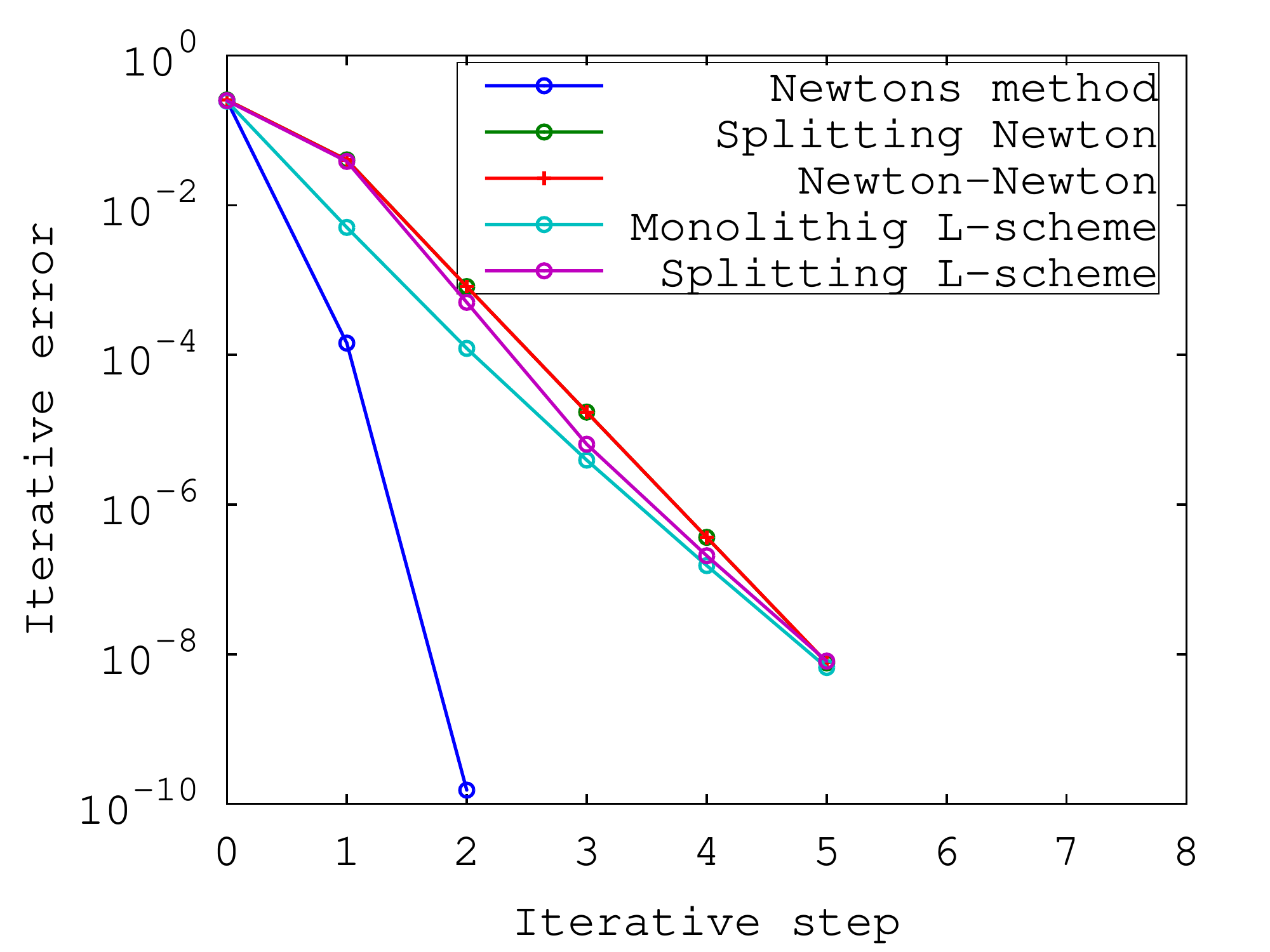} 
\includegraphics[scale=0.25]{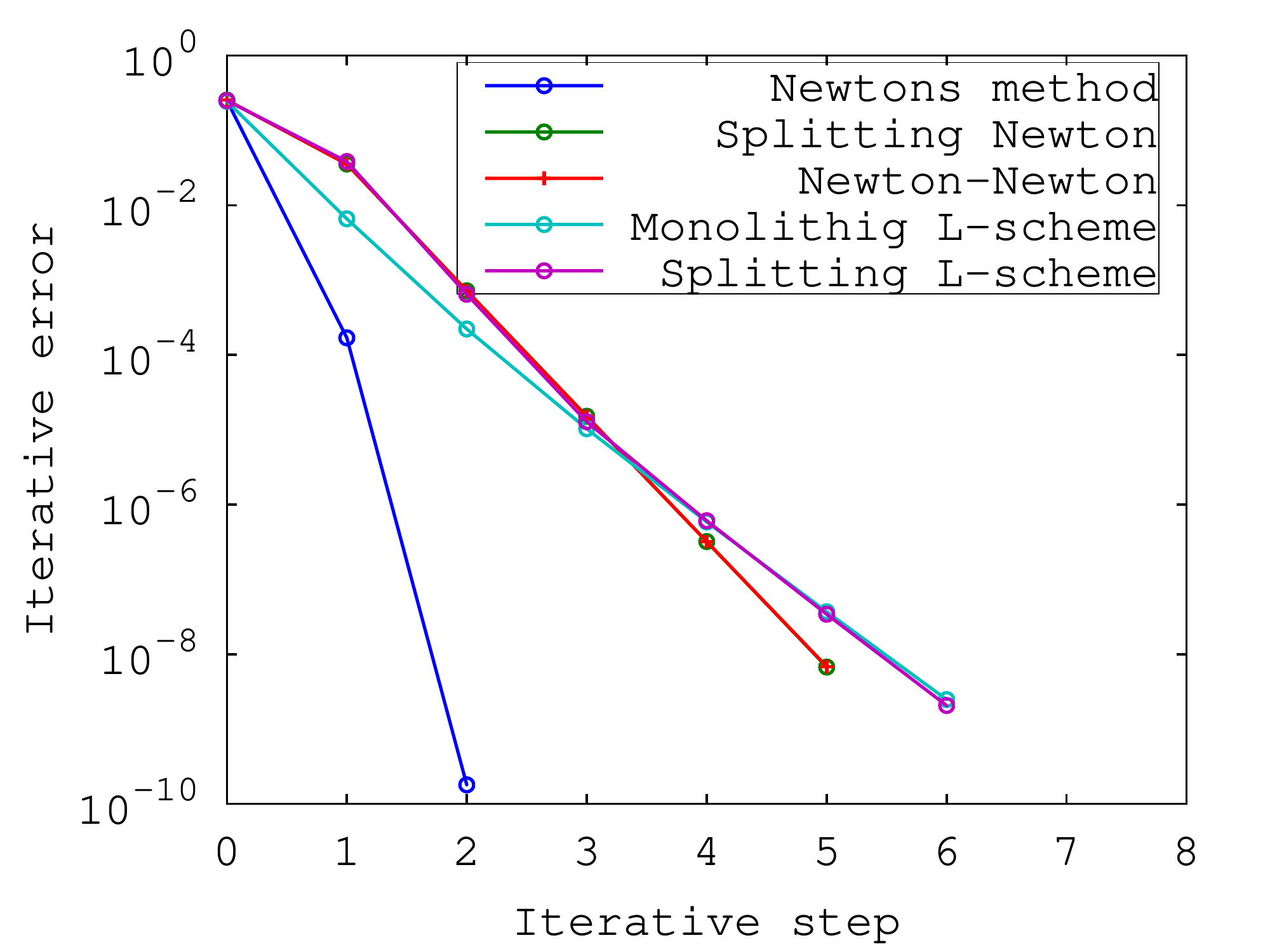} 
\includegraphics[scale=0.25]{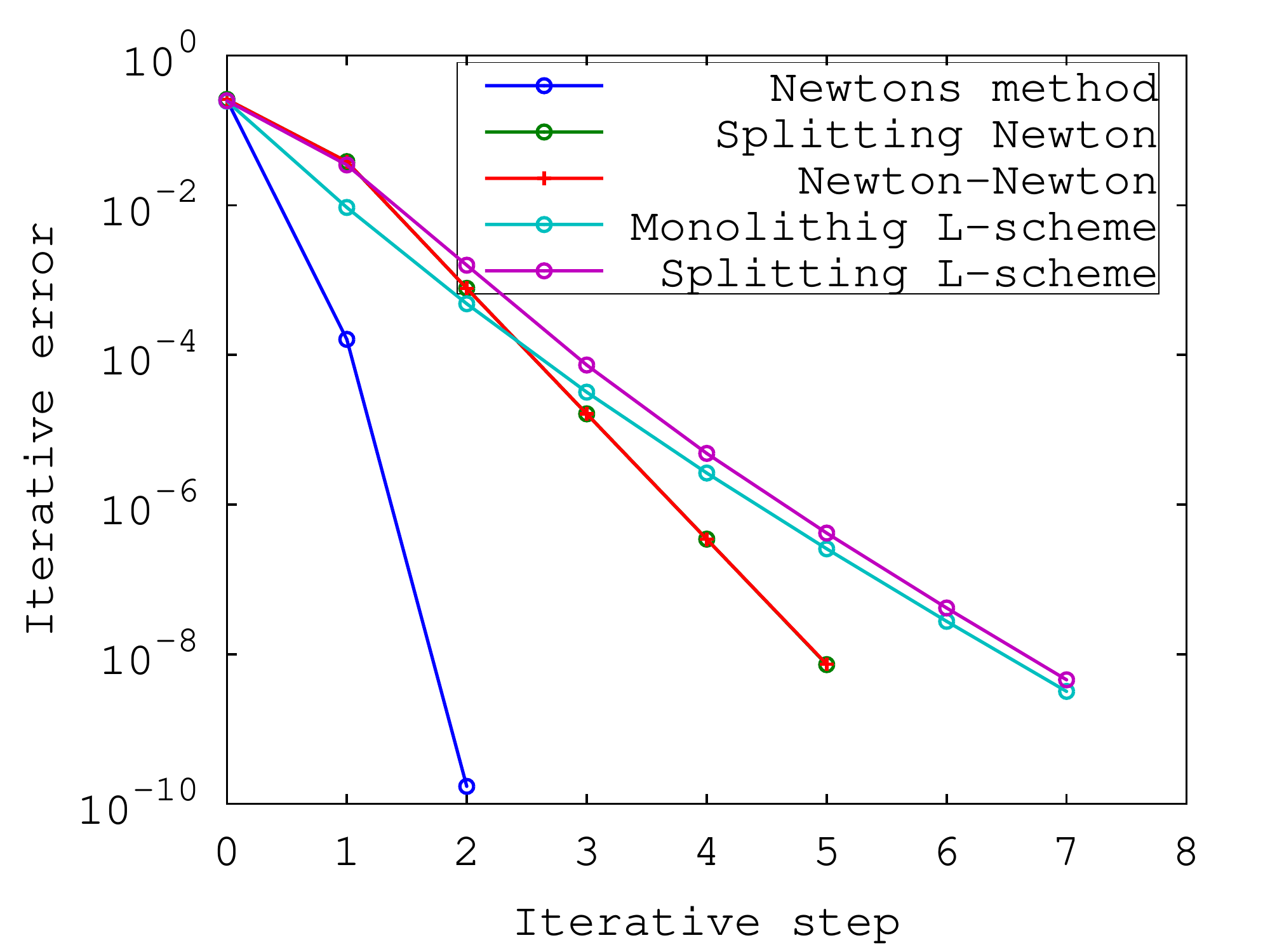} 
\includegraphics[scale=0.25]{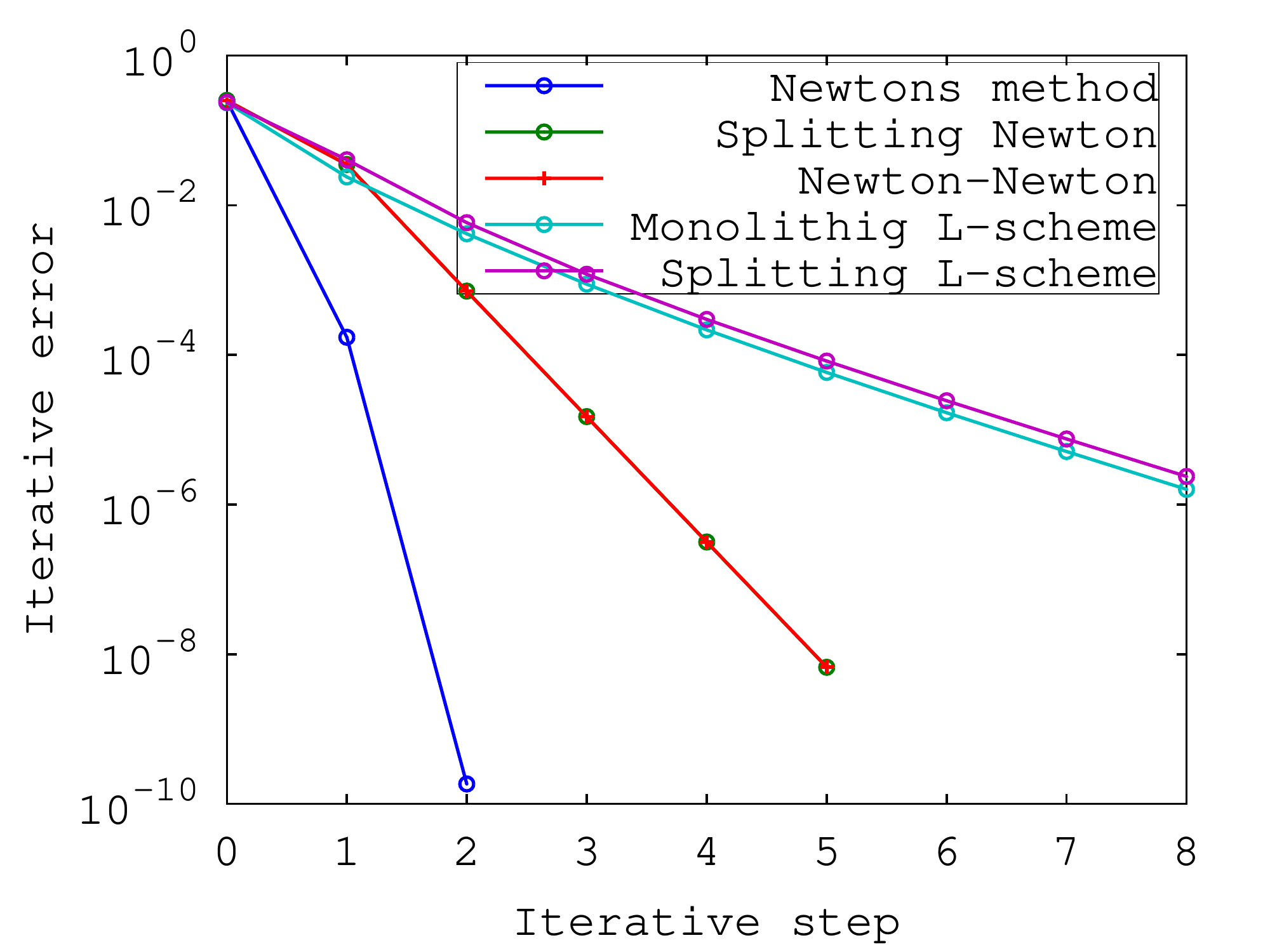} 
\caption{Iterative error at each iteration for different methods:
 to the right $\bb(p) = e^p$, $ \cc(\divv  \uu)  = \sqrt[3]{ \uu^5 }+ \divv \uu$, to the left $\bb(p) = p^2$, $ \cc(\divv  \uu)  = \divv \uu^2 $.}
\label{performances}
\end{figure}

In order to study the performance of the considered schemes, we propose four coefficient functions for $\mathfrak{b(\cdot)}$ and two
for $\mathfrak{c(\cdot)}$, and define four test cases as given in Table \ref{Nonlinearfuntions}. Figure \ref{performances} shows the performance of the numerical methods 
at the last time step $T=1$. The monolithic Newton method shows quadratic convergence in all cases. Nevertheless, the alternate Newton and the $L$-scheme methods show linear convergence as predicted in Section \ref{secc:SmalDef}.

\begin{figure}[h!]
\centering
\includegraphics[scale=0.25]{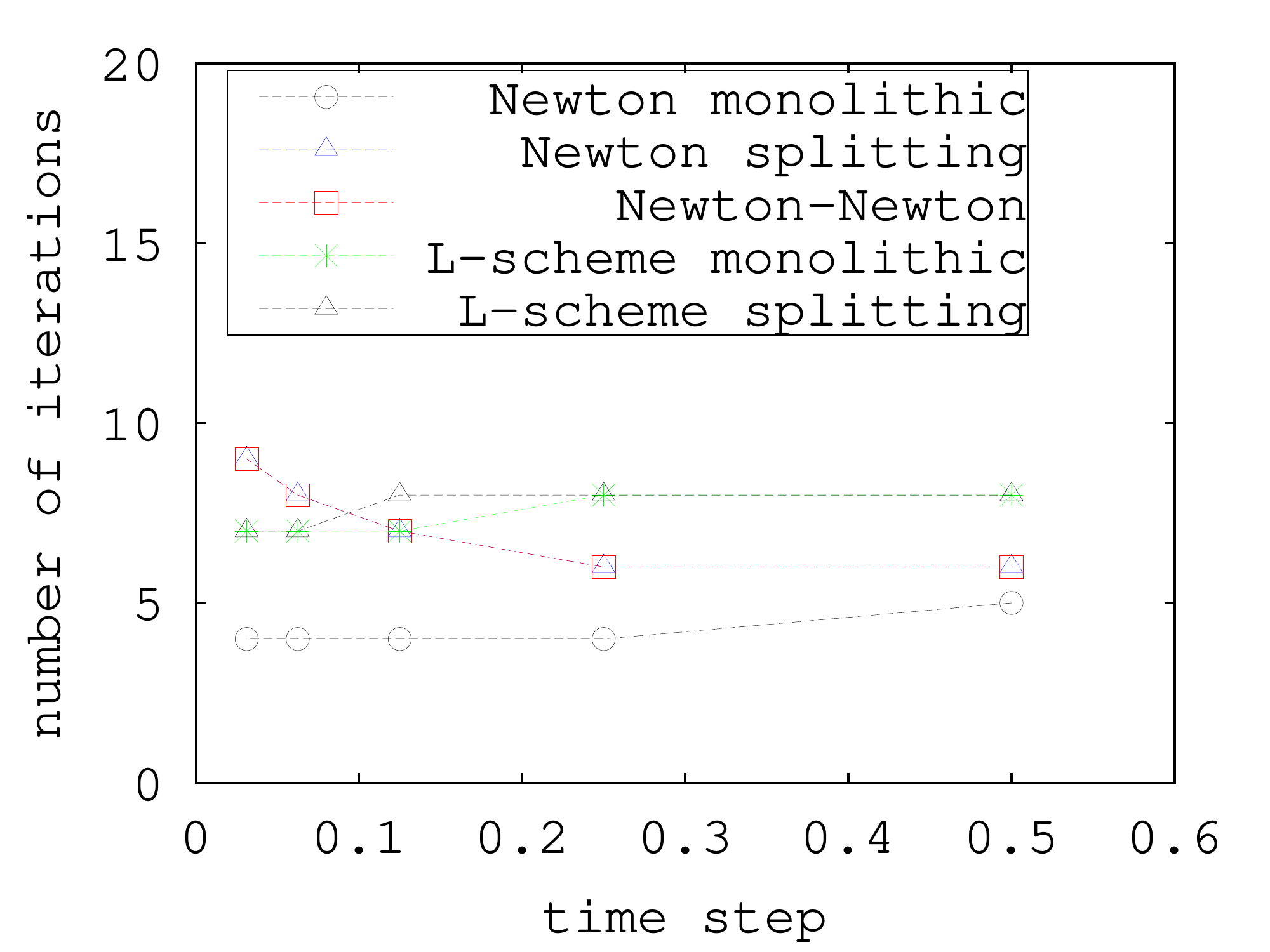} 
\includegraphics[scale=0.25]{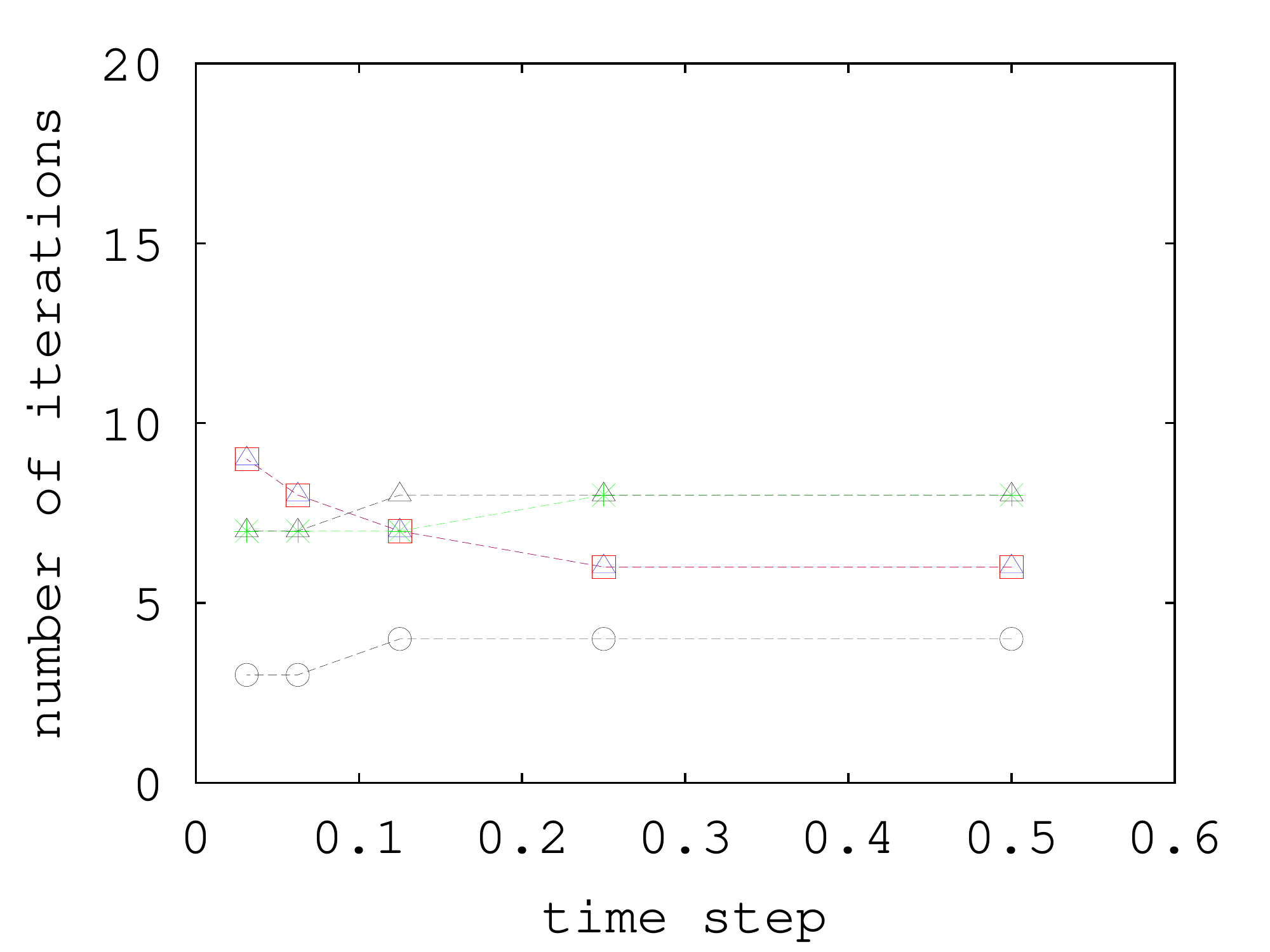} 
\includegraphics[scale=0.25]{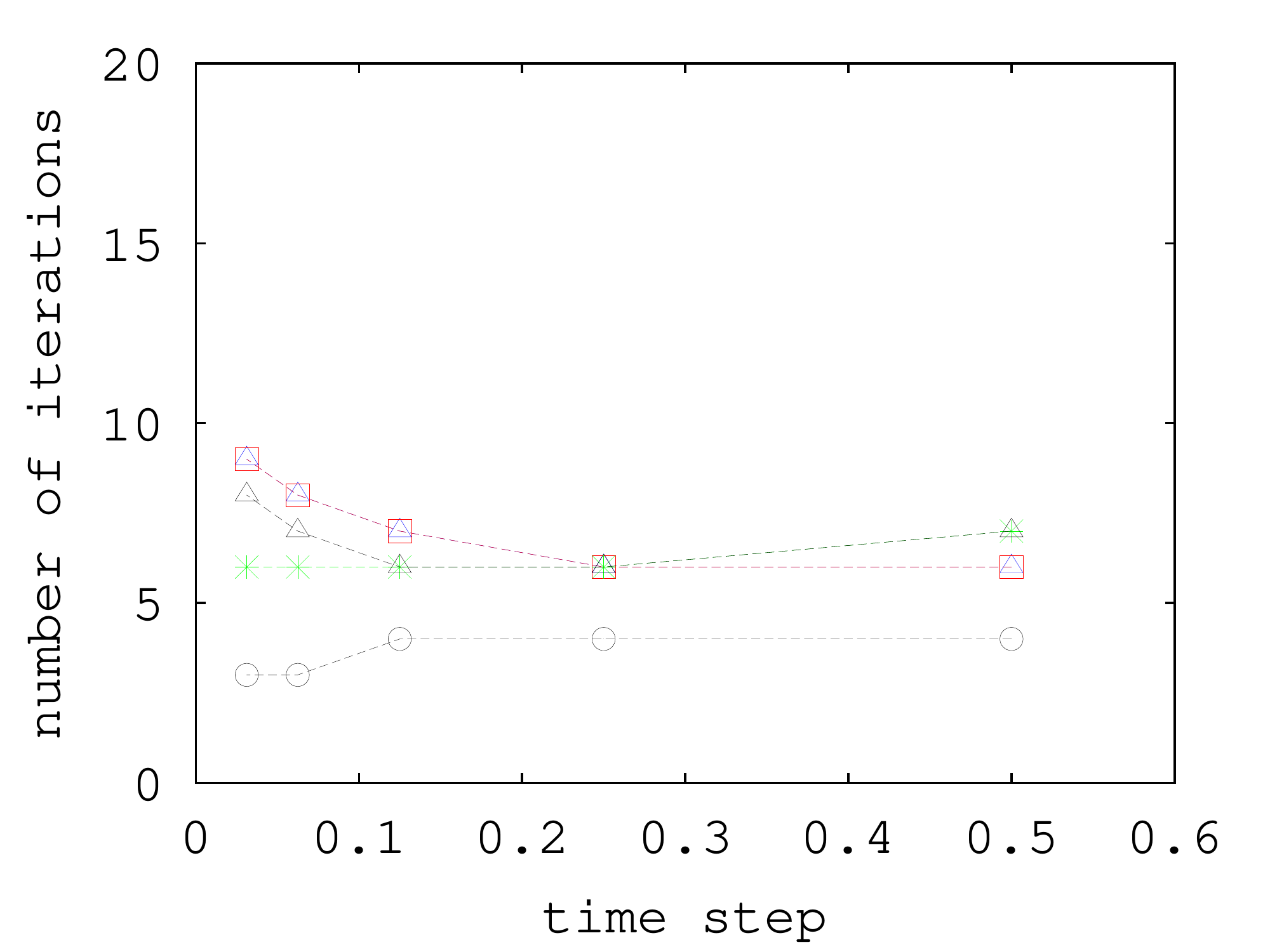} 
\includegraphics[scale=0.25]{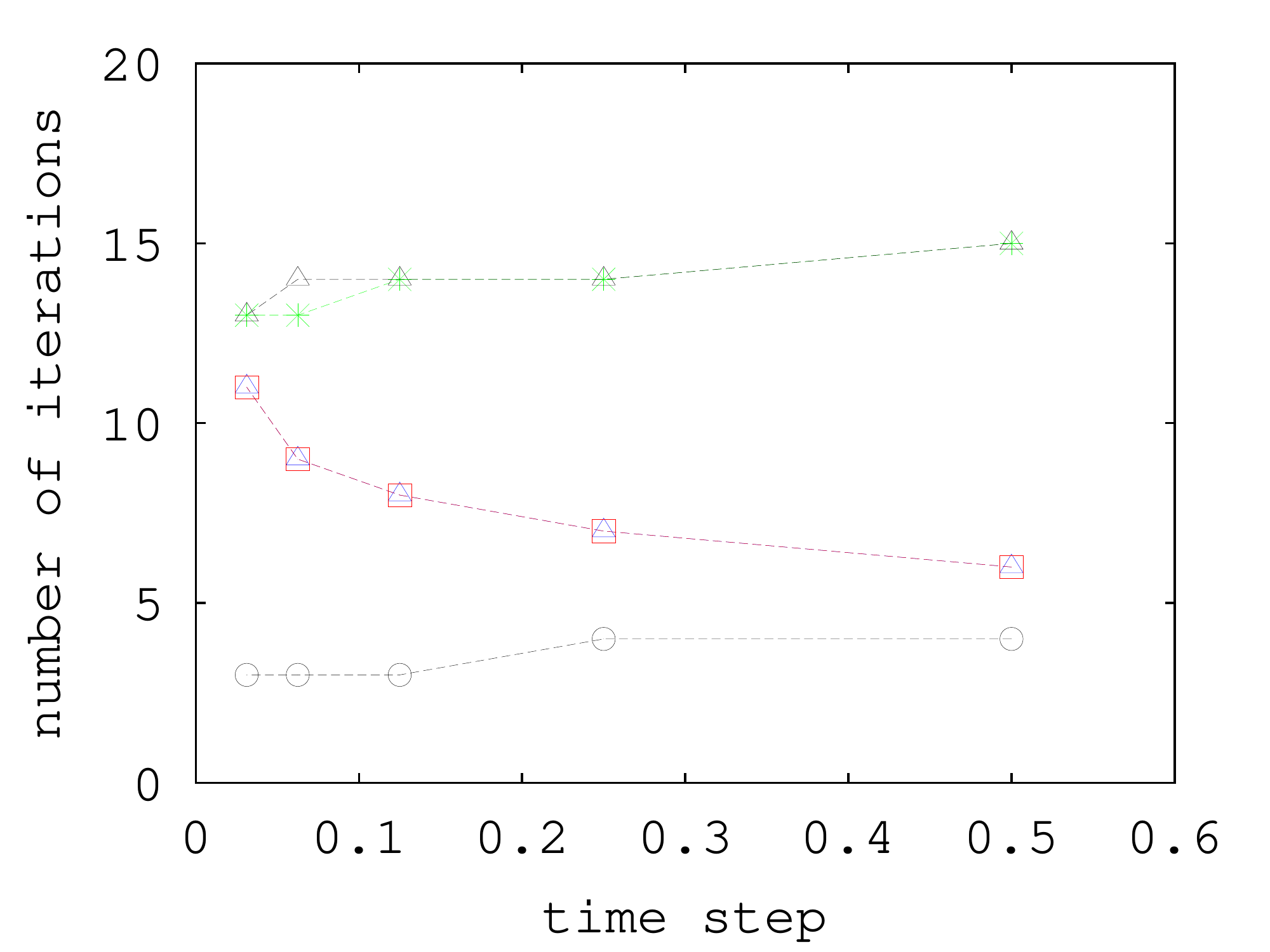} 
\caption{Number of iteration for different time steps:
 to the right $\bb(p) = e^p$, $ \cc(\divv  \uu)  = \sqrt[3]{ \uu^5 }+ \divv \uu$, to the left $\bb(p) = p^2$, $ \cc(\divv  \uu)  = \divv \uu^2 $.}
\label{small_time_step}
\end{figure}

Figure \ref{small_time_step} shows the performance of the considered schemes for different time steps. The Newton method has better convergence for smaller time steps while the $L$-scheme has it for larger time steps; all this is in agreement with the Theorems \ref{Theorem1} and \ref{Newton_Spliting_convergence}.
The performance of the considered schemes are independent of the mesh discretization. 

\begin{table}[h]
 \begin{center}
\caption{The coefficient functions $\bb(\cdot), \cc(\cdot)$ for test problem 1.} 
 \label{Nonlinearfuntions}
 \begin{tabular}{lc|c}
\hline
Case &$\bb(p)$& $\cc(\divv{\uu})$\\
\hline
1  &$e^p$&$ (\divv{\uu})^{3} + \divv{\uu}$\\
2  &$e^p$& $      (\divv{\uu})^{3} $\\
3  &$e^p $&$   \sqrt[3]{\divv{\uu}^{5} } +\divv{\uu}  $\\
4  &$ p^2$&$ \divv{\uu}^{2}$ \\
\hline
\end{tabular}
 \end{center}
\end{table}

\noindent {\bf Test problem 2: a unit cube under large deformation}

We now solve a large deformation problem on the unit-cube  $\Omega=(0,1)^3$. 
A Lagrangian frame of reference is necessary to keep track of the 
deformed domain $\Omega_t$ at time $t$.  We study the performance of the iterative schemes  presented in Section \ref{IterativeSchemes} for solving Eqs. \eqref{darcy}.
The material is supposed to be isotropic and with constant Lam\'e parameters $\mu$ and $\cc(\cdot)$.  We consider a Lagrangian fluid mass  $m_f  = \rho_f J  \phi$ of a slightly compressible fluid, where $\phi$ is the  porosity.  Under this assumption, the time derivative of the fluid content  reads as
 $$\dot{\Gamma}(\uu,p) = c_p J(\uu)\phi \dot{p} +  c_{\alpha} \dot{J}(\uu),  $$ where the compressibility $c_p$ and Biot's coefficient $c_{\alpha} =  J \frac{ \partial \phi}{\partial J} + \phi \approx 1 $ for simplicity. 
We will compare the iterative schemes for a torsion case on a unit cube.
On the top face, we apply the rotation tensor  $\mathbf{R}(\theta)$ of a time dependent angle $\theta (t) =   \pi/4 \ t   $, which gives a rotation of $\pi/4$ at $T=1$. We set homogeneous initial condition for $(\qq_0, p_0)$ and $\nabla \uu_0 = \left(\mathbf{R}(\theta)-\mathbf{I} \right)$. In the alternate Newton method, the stabilization parameter is set to $L_s = 1.$ In the L-scheme method, the linearisation tensor parameters are set as follows: ${\bf{L}}_{\uu} = \mJ_{\uu} {\bf \Pi} \left( \nabla \uu_{0} ,  p_{0} \right),$ ${\bf{L}}_{p} = \mJ_{p} {\bf \Pi} \left( \nabla \uu_{0},  p_{0} \right),$ $\ {\bf{L}}_{\qq} = \mJ_{p} {\bf K} \left( \nabla \uu_{0}\right), $
 $L_p =   \mJ_{p} {\bf \Gamma} \left( \nabla \uu_{0},  p_{0} \right)$ and $L_{\uu} = \mJ_{\uu} {\bf \Gamma} \left( \nabla \uu_{0},  p_{0} \right)$. The mesh size and the time step are set as $h = \tau = 2^{-3}$.  
 We denote by top face of the unit-cube the region $z=1$, the bottom face $ z=0$ and the lateral faces are $x=0$, $x=1$, $y=0$ and $y=1$. The boundary conditions are listed in Table \ref{BoundaryCondition_Traction} and the displacement and pressure field are shown in Figure \ref{Solution_Large_deformation}.

	 \begin{table}[h!]
\centering
\caption{Boundary conditions for Traction and Rotation case respectively.}
  \hspace{1cm}
    \begin{tabular}{ l  l  l }
    \hline
    Face  & Flow & Mechanics \\ 
    Top  & $p=0$ & $\uu = \left(\mathbf{R}(\theta(t))-\bI \right) X_0$ \\
    Bottom &$p=0$& $\uu \cdot  \vec{n}=0$ \\
    Lateral & $p=0$ & $\vec{\Pi} \cdot  \vec{n}=0$ \\
\hline
  \end{tabular}
  \label{BoundaryCondition_Traction}
	\end{table}

\begin{figure}[h!]
\centering
\includegraphics[scale=0.25,trim={0 0 4cm 8cm},clip]{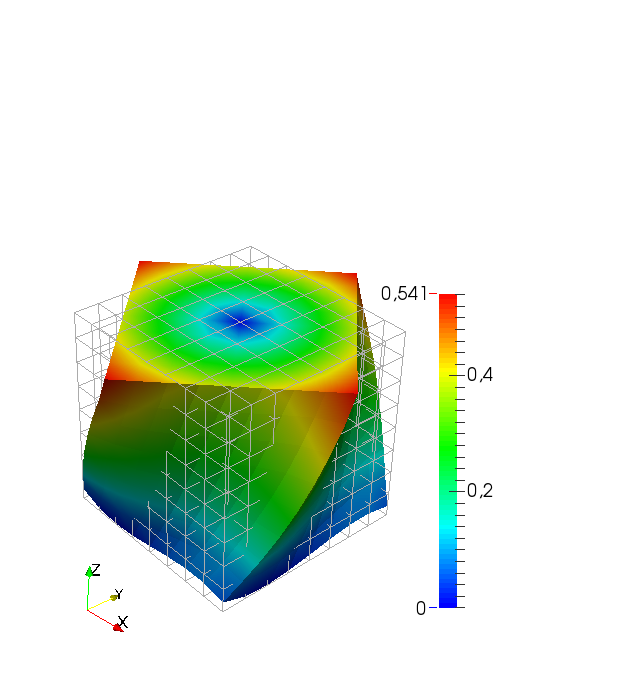} 
\caption{Magnitude of the deformation field and the fluid flow field for torsion.}
\label{Solution_Large_deformation}
\end{figure}	

\begin{figure}[h!]
\centering

\includegraphics[scale=0.30]{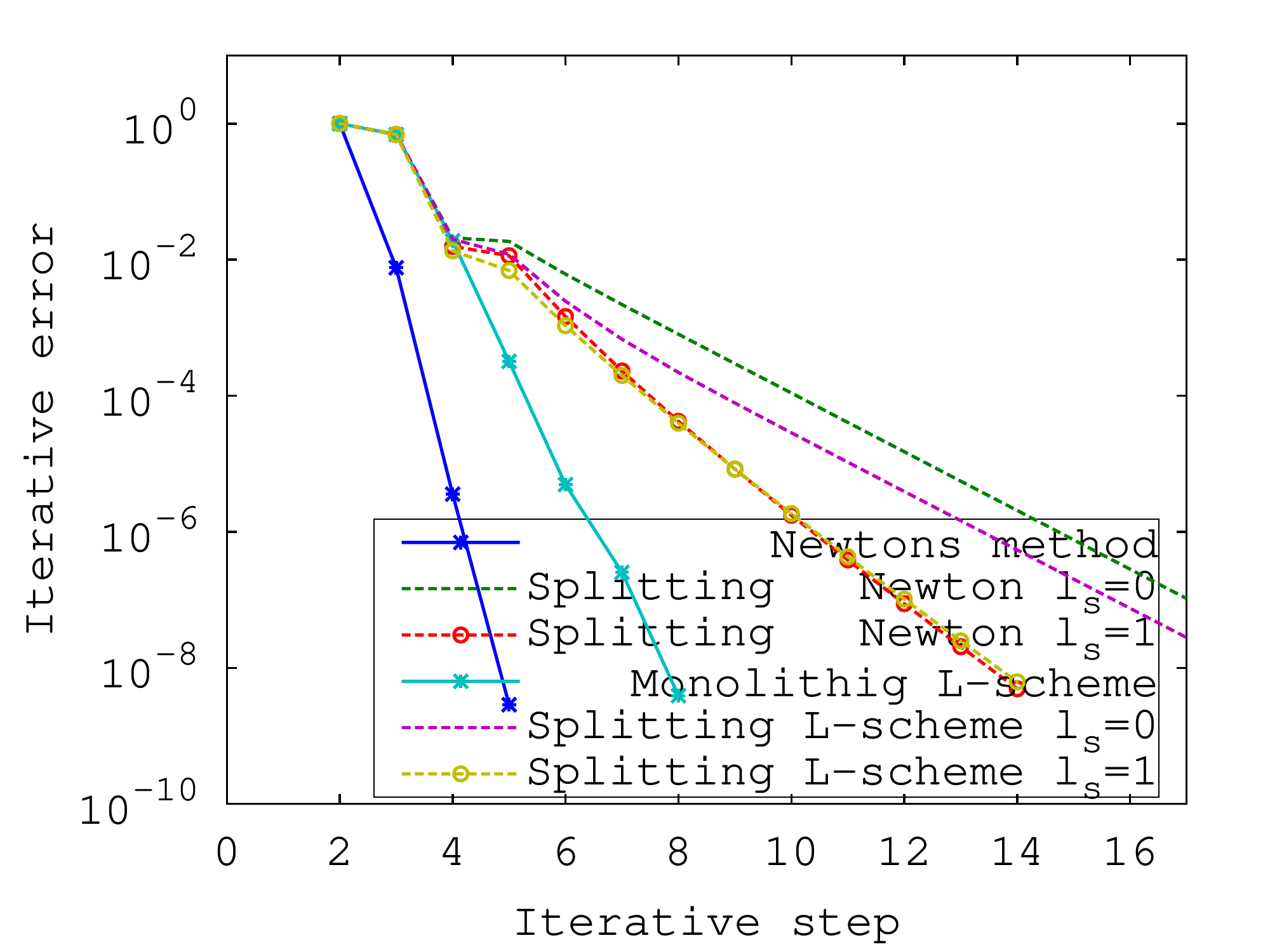} 
\caption{Iterative error at each iteration step for each iterative schemes.}
\label{Solution_Large_deformation}
\end{figure}

We compare the performance of the schemes proposed in Section \ref{IterativeSchemes} and we observe that the numerical convergence is in accordance with the theory developed in Section \ref{secc:SmalDef}, even though the analysis is done for small deformation.
Newton's method has quadratic convergence for the smaller time steps and linear convergence for the larger time steps. In contrast, the monolithic $L$-scheme has the same rate of convergence regardless of size of the time step (see Figure \ref{Solution_Large_deformation}). All splitting schemes have better convergence when the stability term  is used (we use $L_s =1.0$).

\begin{figure}[h!]
\centering
\includegraphics[scale=0.25]{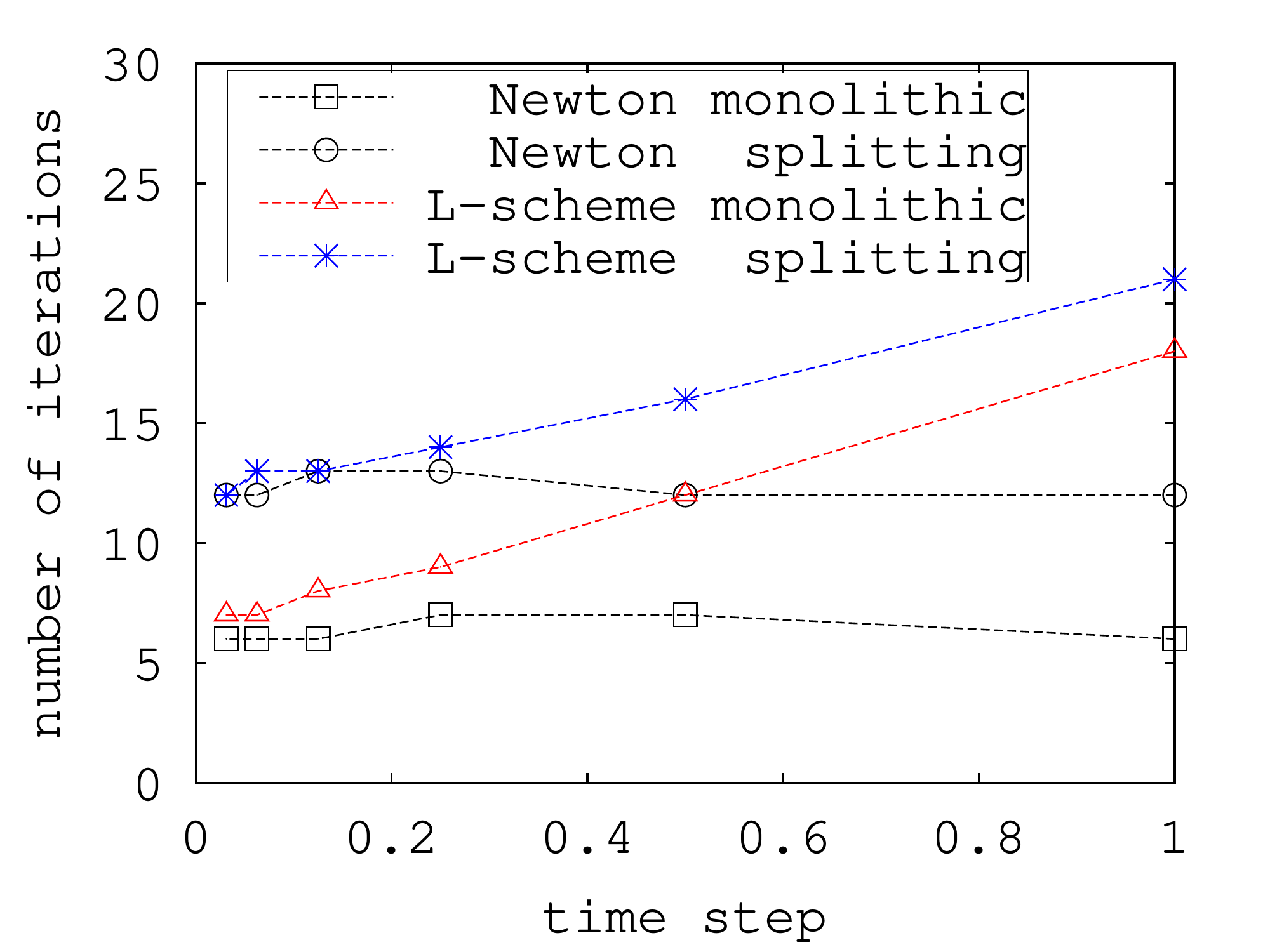} 
\includegraphics[scale=0.25]{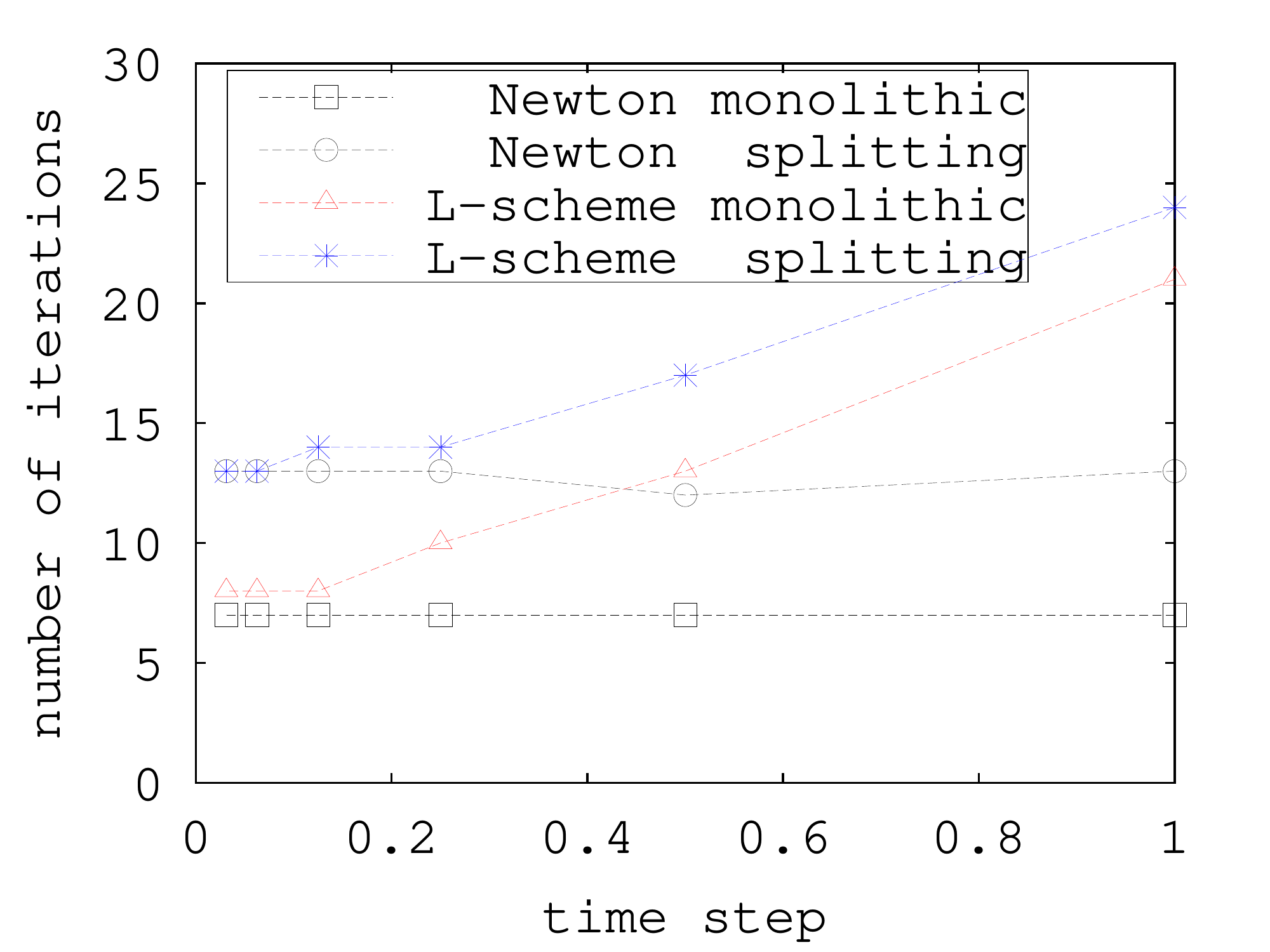} 
\caption{Number of iterations at time $t = 1.0$ using different time steps:
 to the left $h= 1/2^3$, and $h= 1/2^4$ to the right.  }
\label{time_step}
\end{figure}

\section{Conclusions}
\label{Conclusions}

We considered Biot's model under small and large deformation. Different nonlinear solvers based on the $L$-scheme, Newton's method, and the undrained splitting method were presented. The only quadratic convergent scheme is the monolithic Newton method. The splitting Newton method also requires a  stabilization  parameter, otherwise  the  (linear)  convergence  cannot be guaranteed.  The analysis of the schemes and illustrative numerical experiments were presented.

We tested the performance of the schemes on two test problems: a unit square under small deformation and a unit cube under large deformation. To summarise, we make the following remarks:
\begin{itemize}
\item[$\bullet$]{Monolithic and splitting $L$-schemes are robust with respect to the choice of the linearization parameter, the mesh size, and time step size.}
\item[$\bullet$]{The stabilization parameter $L_s$ has a strong influence on the speed of the convergence of the splitting Newton scheme.}
\item[$\bullet$]{The splitting $L$-scheme can be used both as a robust solver or even as a preconditioner (as it is established in \cite{Naga_2010_preconditioner,Tchelepi_2016_Preconditionet}) to improve the performance of the monolithic Newton method and the $L$-scheme.}
\end{itemize}





 	%



\nocite{yi2016}

\clearpage 
\appendix

\section{ Convergence proof of the alternate Newton method}
\label{Appendix_A} 
 
The following result provides the linear convergence of the alternate Newton method in~\eqref{small_splitting2_Newton}-~\eqref{small_splitting1_Newton} for $\tau$ sufficiently small.

\begin{theorem}
Assuming (A1)-(A4) and $L_s\geq \frac{\alpha^2}{\alpha_{\bb}}$, the alternate Newton splitting method in~\eqref{small_splitting2_Newton}-\eqref{small_splitting1_Newton} converges linearly if $\tau$ is small enough.
\end{theorem} 
 \proof
By subtracting problems~\eqref{small_splitting2_Newton}-\eqref{small_splitting1_Newton} and  \eqref{fully1_discret}, taking as test functions  ${\mathbf e}_{\qq}^{n,i}$, $e_p^{n,i}$ and ${\mathbf e}_{\mathbf u}^{n,i}$, and rearranging some elements to the right hand side we obtain,

\begin{align}
\begin{array}{rl}
\left(\bK^{-1}{\bf e}_{\qq}^{n,\ii} , {\bf e}_{\qq}^{n,\ii} \right) - \left( e_p^{n,\ii} , \divv \be_{\qq}^{n,\ii}\right)&= 0,
\label{splitting3_error} 
\end{array}
\\
\begin{array}{rl}
\label{splitting33_error} 
 \left( \bb'(p_h^{n,i-1})(p_h^n-p_h^{n,i}),e_p^{n,i} \right) + \alpha \left(   \displaystyle \nabla \cdot {\mathbf e}_{\mathbf u}^{n,i-1}, e_p^{n,i} \right)  
+ \tau \left(\divv \be_{\qq}^{n,\ii}, e_p^{n,\ii} \right)
 \\ = \left( \bb(p_h^n) - \bb(p_h^{n,i-1})-\bb'(p_h^{n,i-1})(p_h^{n}-p_h^{n,i-1}), e_p^{n,i} \right).
\end{array}
\end{align}

The mechanics equation  then gives,
\begin{align}
 \label{splitting2_error} 
\begin{array}{rl}
 \left(\varepsilon ({\mathbf e}_{\mathbf u}^{n,i}), \varepsilon ({\mathbf e}_{\mathbf u}^{n,i})\right) +
\left(\cc'(\nabla \cdot {\mathbf u}_h^{n,i-1}) 
\nabla \cdot  {\mathbf e}_{\mathbf u}^{n,i} ,\nabla \cdot {\mathbf e}_{\mathbf u}^{n,i} \right)  
    \\ [1.3ex]
  + L_s \left( \nabla \cdot \delta {\bu}_h^{n,i}, \nabla\cdot {\mathbf e}_{\mathbf u}^{n,i} \right)  
- \alpha \left(  e_p^{n,i}, \nabla\cdot {\mathbf e}_{\mathbf u}^{n,i}\right) 
  \\ [1.3ex]
 = \left( \cc(\nabla \cdot {\mathbf u}_h^n)-\cc(\nabla \cdot {\mathbf u}_h^{n,i-1}) +  \cc'(\nabla \cdot {\mathbf u}_h^{n,i-1}) 
\nabla \cdot {\mathbf e}_{\mathbf u}^{n,i-1},\nabla \cdot {\mathbf e}_{\mathbf u}^{n,i} \right).  
\end{array}
\end{align}
By using similar steps as in Theorem  \ref{Theorem1}, we obtain the following

\begin{align}
\label{splitting_first_bound_g}
\begin{array}{rl}
||\varepsilon({\mathbf e}_{\mathbf u}^{n,i})||^2+\left(\cc'(\nabla\cdot {\mathbf u}_h^{n,i-1})
\nabla \cdot {\mathbf e}_{\mathbf u}^{n,i},\nabla \cdot {\mathbf e}_{\mathbf u}^{n,i}\right)
 \\ [1.3ex]
+ L_s \left(\nabla \cdot ({\mathbf e}_{\mathbf u}^{n,i}- {\mathbf e}_{\mathbf u}^{n,i-1}) , \nabla\cdot {\mathbf e}_{\mathbf u}^{n,i}\right)
-\alpha \left( e_p^{n,i},\nabla \cdot {\mathbf e}_{\mathbf u}^{n,i}\right)  
 \\ [1.3ex]
\leq \displaystyle \frac{L_{\cc'}^2}{8\gamma_1}||\nabla \cdot {\mathbf e}_{\mathbf u}^{n,i-1}||_{L^4(\Omega)}^4 + \frac{\gamma_1}{2}||\nabla \cdot {\mathbf e}_{\mathbf u}^{n,i}||^2. 
\end{array}
\end{align}
Next, by using the inverse inequality $||\cdot||_{L^4(\Omega)}\leq C h^{-d/4}||\cdot||$ \cite{Brenner_1991}, and by using the following formula  $(x-y,x) = \displaystyle \frac{||x||^2}{2} +\frac{||x-y||^2}{2}-\frac{ ||y||^2}{2}$,  by choosing $x = \nabla\cdot {{\mathbf e}_{\mathbf u}}^{n,i}$ and $y = \nabla\cdot {\mathbf e}_{\mathbf u}^{n,i-1},$ 
  we obtain from~\eqref{splitting_first_bound_g}

\begin{align}
\label{splitting_3_bound_g}
\begin{array}{rl}
||\varepsilon({\mathbf e}_{\mathbf u}^{n,i})||^2+\left(\cc'(\nabla\cdot {\mathbf u}_h^{n,i-1})
\nabla \cdot {\mathbf e}_{\mathbf u}^{n,i},\nabla \cdot {\mathbf e}_{\mathbf u}^{n,i}\right)  + \frac{L_s}{2} \norm{\nabla \cdot ({\mathbf e}_{\mathbf u}^{n,i}- {\mathbf e}_{\mathbf u}^{n,i-1}) }^2 
 \\ [1.3ex]
\frac{L_s}{2} \norm{ \nabla\cdot {\mathbf e}_{\mathbf u}^{n,i}}^2
- \alpha\left( e_p^{n,i},\nabla \cdot {\mathbf e}_{\mathbf u}^{n,i}\right)+
\leq \displaystyle C_1 h^{-d}\frac{L_{\cc'}^2}{8\gamma_1}||\nabla \cdot {\mathbf e}_{\mathbf u}^{n,i-1}||^4 
 \\ [1.3ex]
+ \frac{\gamma_1}{2}||\nabla \cdot {\mathbf e}_{\mathbf u}^{n,i}||^2 
+\frac{L_s}{2}\norm{\nabla\cdot {\mathbf e}_{\mathbf u}^{n,i-1}}^2.
\end{array}
\end{align}
Finally,  by reorganizing \eqref{splitting_3_bound_g},
using (A2) and choosing $\gamma_1=\alpha_{\cc}$, we obtain the following inequality,

\begin{align} \label{splitting_5_bound_g}
\begin{array}{rl}
||\varepsilon({\mathbf e}_{\mathbf u}^{n,i})||^2+
\left( \frac{\alpha_{\cc}+ L_s}{2}\right) \norm{ \nabla\cdot {\mathbf e}_{\mathbf u}^{n,i}}^2 + \frac{L_s}{2} \norm{\nabla \cdot \delta {\mathbf e}_{\mathbf u}^{n,\ii} }^2
\\ [1.3ex]
\leq \displaystyle C_1 h^{-d}\frac{L_{\cc'}^2}{8\alpha_{\cc}}||\nabla \cdot {\mathbf e}_{\mathbf u}^{n,i-1}||^4 
 +\frac{L_s}{2}\norm{\nabla\cdot {\mathbf e}_{\mathbf u}^{n,i-1}}+ \alpha\left( e_p^{n,i},\nabla \cdot {\mathbf e}_{\mathbf u}^{n,i}\right).
 \end{array}
\end{align}
In a similar way, we obtain the following expression from~\eqref{fully2_error},
\begin{equation}\label{splitting_final_bound_b}
\frac{\tau}{k_M}||\be_{\qq}^{n,i}||^2+\frac{\alpha_{\bb}}{2}||e_p^{n,i}||^2
 \leq \displaystyle C_2 h^{-d}\frac{L_{\bb'}^2}{8\alpha_{\bb}}||e_p^{n,i-1}||^4 -\alpha \left(\nabla \cdot {\mathbf e}_{\mathbf u}^{n,i-1}, e_p^{n,i}\right).
\end{equation}
Adding equations~\eqref{splitting_5_bound_g} and~\eqref{splitting_final_bound_b} yields,
\begin{align}
\label{splitting_1_adding_bound}
\begin{array}{rl}
\frac{\tau}{k_M}||\be_{\qq}^{n,i}||^2+\frac{\alpha_{\bb}}{2}||e_p^{n,i}||^2 
+||\varepsilon({\mathbf e}_{\mathbf u}^{n,i})||^2
+ \frac{L_s}{2} \norm{\nabla \cdot \delta {\mathbf e}_{\mathbf u}^{n,i} }^2
+
\left( \frac{\alpha_{\cc}+ L_s}{2}\right) \norm{ \nabla\cdot {\mathbf e}_{\mathbf u}^{n,i}} ^2
 \\ [1.3ex]
\leq \displaystyle C_2 h^{-d}\frac{L_{\bb'}^2}{8\alpha_{\bb}}||e_p^{n,i-1}||^4+ \displaystyle C_1 h^{-d}\frac{L_{\cc'}^2}{8\alpha_{\cc}}||\nabla \cdot {\mathbf e}_{\mathbf u}^{n,i-1}||^4  
 \\[1.3ex]
+\frac{L_s}{2}\norm{\nabla\cdot {\mathbf e}_{\mathbf u}^{n,i-1}}^2
+\alpha\left( \nabla \cdot \delta  {\mathbf e}_{\mathbf u}^{n,i} ,e_p^{n,i}\right).
\end{array}
\end{align}

By using Young's inequality $(a,b)\leq \displaystyle \frac{||a||^2}{2\gamma}+\frac{\gamma ||b||^2}{2}$, for $\gamma > 0$ and  choosing $b = e_p^{n,i}$ and $a = \nabla \cdot \delta  {\mathbf e}_{\mathbf u}^{n,i} $ we bound the coupling term (for $\gamma_2>0$), 

\begin{equation}
\alpha\left( \nabla \cdot  \delta  {\mathbf e}_{\mathbf u}^{n,i} ,e_p^{n,i}\right) \leq   \frac{\alpha^2}{2\gamma_2} \norm{ \nabla \cdot \delta  {\mathbf e}_{\mathbf u}^{n,i-1} }^2 +\frac{ \gamma_2}{2}\norm{e_p^{n,i}}^2.
\label{intermediate_coplinterm}
\end{equation}
Then by using \eqref{intermediate_coplinterm} and choosing $\gamma_2 = \frac{\alpha_{\bb}}{2}$ we obtain from \eqref{splitting_1_adding_bound}

\begin{align}
\label{splitting_final_bound_1}
\begin{array}{rl}
\frac{\tau}{k_M}||\be_{\qq}^{n,i}||^2+  \frac{\alpha_{\bb}}{4}||e_p^{n,i}||^2 
+||\varepsilon({\mathbf e}_{\mathbf u}^{n,i})||^2
+ \left(  \frac{L_s}{2}-\frac{\alpha^2}{2\alpha_{\bb}} \right) \norm{\nabla \cdot \delta {\mathbf e}_{\mathbf u}^{n,i} } ^2
 \\[1.3ex]
+
\left( \frac{\alpha_{\cc}+ L_s}{2}\right) \norm{ \nabla\cdot {\mathbf e}_{\mathbf u}^{n,i}} ^2
\leq   \frac{h^{-d}}{8}\left( C_2\frac{L_{\bb'}^2}{\alpha_{\bb}}||e_p^{n,i-1}||^4+ C_1 \frac{L_{\cc'}^2}{\alpha_{\cc}}||\nabla \cdot {\mathbf e}_{\mathbf u}^{n,i-1}||^4 \right)
 \\[1.3ex]
 +\frac{L_s}{2}\norm{\nabla\cdot {\mathbf e}_{\mathbf u}^{n,i-1}}^2. 
\end{array}
\end{align}
Since  $L_s\geq \frac{\alpha^2}{\alpha_{\bb}}$, we obtain  
\begin{align}
\label{Asplitting_final_bound}
\begin{array}{rl}
 \frac{\tau}{k_M}||\be_{\qq}^{n,i}||^2+  \frac{\alpha_{\bb}}{4}||e_p^{n,i}||^2 
+
\left( \frac{\alpha_{\cc}+ L_s}{2}\right) \norm{ \nabla\cdot {\mathbf e}_{\mathbf u}^{n,i}} ^2  
\\ [1.3ex]
\leq   \frac{h^{-d}}{8}\left(C_2 \frac{L_{\bb'}^2}{\alpha_{\bb}}||e_p^{n,i-1}||^4+ C_1\frac{L_{\cc'}^2}{\alpha_{\cc}}||\nabla \cdot {\mathbf e}_{\mathbf u}^{n,i-1}||^4 \right)
 \\ [1.3ex]
 +\frac{L_s}{2}\norm{\nabla\cdot {\mathbf e}_{\mathbf u}^{n,i-1}}^2. 
\end{array}
\end{align}
\noindent By using $\norm{\nabla \cdot {\mathbf e}_{\mathbf u}^{n,0}}\leq C\tau$, $\norm{ {e}_{p}^{n,0}}\leq C  \tau$ wich can be proven and the estimate in  Lemma \ref{succeion_convergence}, the convergence is ensured if  $\tau = O (h^\frac{d}{2})$.





\endproof

\section*{Acknowledgement}
The research was supported by the University of Bergen in cooperation with the FME-SUCCESS center (grant 193825/S60) funded by the Research Council of Norway.  The work has also been partly supported by: the NFR-DAADppp grant 255715, the NFR-Toppforsk project 250223, the NRC-CHI grant  255510 and the NRC-IMMENS grant 255426. 
 
 \bibliographystyle{plain}

\end{document}